\documentclass[11pt,reqno]{amsart}
\usepackage{amsmath}
\usepackage{amssymb}
\usepackage[sans]{dsfont}

\def\squarebox#1{\hbox to #1{\hfill\vbox to #1{\vfill}}}

\theoremstyle{plain}

\newtheorem{thm}{Theorem}
   
\newtheorem{cor}{Corollary}
   
\newtheorem{lem}{Lemma}

\DeclareMathOperator{\re}{Re}
\DeclareMathOperator{\im}{Im}

\renewcommand{\Re}{\mathop{\rm Re}\nolimits}
\renewcommand{\Im}{\mathop{\rm Im}\nolimits}

\def\la{\langle}
\def\ra{\rangle}

\def\one{\mathds{1}}
\def\Oc{{\mathcal O}}

\newtheorem{rem}{Remark}
\newtheorem{prop}{Proposition}

\numberwithin{equation}{section}



\begin{document}

\def\C{{\mathbb C}}
\def\R{{\mathbb R}}
\def\N{{\mathbb N}}
\def\Z{{\mathbb Z}}
\def\T{{\mathbb T}}
\def\Q{{\mathbb Q}}
\def\SP{{\mathbb S}}
\def\d{{\partial}}
\def\mc{{\mathcal H}}
\def\1b{{\mathbb I}}
\def\tr{{\rm tr}\:}
\def\pv{\partial_x V}
\def\pa{\partial}
\def\pc{{\mathcal P}}
\def\fm{\frac{1}{1+m}}
\def\om{\omega_0}
\def\h{{\mathcal H}}
\def\p{{\mathcal P}}
\def\ep{\epsilon}

\title[Upper bound]
{Upper bound for the counting function of interior transmission eigenvalues}
\date{}
\author[M. Dimassi]{Mouez Dimassi}
\author[V. Petkov]{Vesselin Petkov}

\address {Universit\'e Bordeaux I, Institut de Math\'ematiques de Bordeaux,  351, Cours de la Lib\'eration, 33405  Talence, France}
\email{Mouez.Dimassi@math.u-bordeaux1.fr}
\address {Universit\'e Bordeaux I, Institut de Math\'ematiques de Bordeaux,  351, Cours de la Lib\'eration, 33405  Talence, France}
\email{petkov@math.u-bordeaux1.fr}
\thanks{The second author was partially supported by the ANR project NOSEVOL}
\maketitle
\begin{abstract} For the complex interior transmission eigenvalues (ITE) we study for small $\theta > 0$ the counting function
$$N(\theta, r) = \#\{\lambda \in \C:\: \lambda \: {\rm is} \: {\rm (ITE)},\: |\lambda| \leq r, \:  |\arg \lambda| \leq \theta\}.$$
We obtain for fixed $\theta > 0$ an upper bound $N(\theta, r) \leq C r^{n/2}, \: r \geq r(\theta).$
\end{abstract}

\date{today}

\section{Introduction}
\renewcommand{\theequation}{\arabic{section}.\arabic{equation}}
\setcounter{equation}{0}

Let $\Omega \subset \R^n,\: n \geq 2$ be a bounded connected open domain with smooth boundary $\pa \Omega.$ A complex number $\lambda \neq 0$ is called an interior  transmission eigenvalue (ITE) if and only if there exist $u \neq 0,\: v \neq 0$, $ u,  v \in H^2(\Omega),$ so that
\begin{equation} \label{eq:1.1}
\begin{cases}  - \Delta u - \lambda u = 0, x \in \Omega,\\
-\Delta v - \lambda(1 + m(x)) v = 0,\: x \in \Omega,\\
\gamma_0 (u - v) = 0, \:\gamma_0\pa_{\nu} (u - v) = 0. \end{cases}
\end{equation}
Here $\nu$ is the exterior unit normal to $\partial \Omega$ and $\gamma_0$ is the operator of restriction on $\partial \Omega.$
In this work we suppose that $m(x) \in C^{\infty}(\bar{\Omega})$,\ is a real-valued function satisfying the conditions
\begin{equation} \label{eq:1.2}
1 + m(x) \geq \delta_0 > 0,\: \forall x \in \bar{\Omega},\: m(x) \neq 0, \: \forall x \in \pa \Omega.
\end{equation}

Consider in the space $\h = L^2(\Omega) \oplus L^2(\Omega)$ the operator
$$\p = \begin{pmatrix} -\Delta & & 0\\ 0 & & -\fm\Delta\end{pmatrix}$$
with domain
$$D(\p) = \{ (u, v) \in \h:\: \Delta u \in L^2(\Omega), \: \Delta v \in L^2(\Omega),\: (u- v) \big\vert_{\partial \Omega} = 0,\:\partial_{\nu} (u - v) \big \vert_{\partial \Omega} = 0\}.$$
This operator is closed (see Lemma 1) and the (ITE) are the eigenvalues of $\p$ in $\C \setminus \{0\}.$\\

 The (ITE) are related to the scattering theory. In particular, if $\lambda = k^2 \in \R$ is not a real (ITE), then the far-field operator $F(\lambda): L^2({\mathbb S}^{n-1}) \longrightarrow L^2({\mathbb S}^{n-1})$ with kernel the scattering amplitude $a(\lambda, \theta, \omega)$ is injective and has a dense range. The (ITE) play an important role in the inverse scattering methods and there is an increasing interest on the problems of the existence and distribution of (ITE). The reader may consult the survey \cite{CH} for a comprehensive review of recent results and references.\\
 
There are several works proving that under some conditions the (ITE) form a discrete set in $\C$ (see \cite{CGH}, \cite{HK1}, \cite{HK2}, \cite{CH}, \cite{Sy}). In the case when $m(x)$ satisfies (\ref{eq:1.2}) the results in \cite{Sy}, \cite{R} guarantee that (ITE) are discrete with finite dimensional generalized eigenspaces. There are many works on the distribution of the (ITE)  (see \cite{BP}, \cite{LV1}, \cite{LV2}, \cite{LV3}, \cite{SS}, \cite{HK3}, \cite{R} and the references cited there). The distribution of the positive (ITE) has been studied in \cite{LV2}, \cite{LV3}, \cite{SS} and a lower bound for the number of the positive (ITE) has been established. The analysis of the complex (ITE)  is more difficult than that of the positive ones. Under the assumption (1.2), it has been proved in \cite{V} that for every $\ep > 0$ outside a parabolic neighborhood of $\R^+$ having the form
$$\{\lambda \in \C:\: |\im \lambda| \leq C_{\ep} (|\re \lambda| + 1)^{3/4 + \ep},\: \re \lambda > 0\}$$
there are at most a finite number of (ITE) (see also \cite{HK3} for an weaker result). From this it follows that the (ITE) can have accumulation point only at infinity and for every fixed $0 < \theta < \pi/2$ outside the domain
$$D(\theta) = \{\lambda \in \C:\:|\arg \lambda| \leq \theta\}$$
 we have at most finite number (ITE) (see also \cite{LV4} for another proof of this result). The multiplicity of an (ITE) $\lambda$ is given by the dimension of the eigenspace $E_{\lambda}$
generated by the eigenfunctions and the generalized eigenfunctions corresponding to $\lambda.$ In this direction we mention the interesting work \cite{HS}, where the case of a constant function $m(x) = m$ and the ball $\{ x \in \R^n: |x| \leq 1\}$ has been studied. This work shows that even in the case $n = 1$ if $m = p/ q$ is rational and $p/q$ is irreducible, the geometric and algebraic multiplicities of an (ITE) may be different (see Theorem 2.4 and 2.5 in \cite{HS}). The reader should consult also  \cite{Sy1} for the analysis of the case when $\Omega = \{|x| \leq 1\}$ and $n=1.$ In the following in the definition of the counting function
$N(r) = \#\{ \lambda \in \C:\: \lambda\:\:{\rm is}\:\: (ITE),\: |\lambda| \leq r\}$
every eigenvalue $\lambda$ will be counted with its {\it algebraic} multiplicity. Moreover, it follows from  Theorem 3, \cite{R}, that under the assumptions (\ref{eq:1.2}) the set of the eigenvalues of $\p$ is infinite and it is natural to study the asymptotics of $N(r)$ as $r \to \infty$. \\ 

In \cite{LV1} a Weyl type asymptotics has been obtained in the {\it anisotropic case}, when in (\ref{eq:1.1}) for $v$ we have the equation $-\nabla A(x) \nabla v - \lambda (1 + m(x)) v = 0$  with a positive definite matrix $A(x)$ different from the identity. Under some conditions on $A(x)$, which guarantee that the boundary problem is parameter-elliptic (see \cite{AV}), it is possible  to apply the results in \cite{BK1}. More precisely, in \cite{LV1} it is shown that the counting function $N(r)$ has the asymptotics
\begin{equation} \label{eq:1.3}
N(r) \sim \frac{\omega_n}{(2\pi)^n} \int_{\Omega} \Bigl(1 + \frac{(1 + m(x))^{n/2}}{(\det A(x))^{1/2}}\Bigr) dx\, r^{n/2},\: r \to \infty,
\end{equation}
where $\omega_n$ is the volume of the unit ball in $\R^n.$ It is important to note that in the case treated in \cite{LV1} the estimates for the solution in \cite{AV} imply that the domain of the corresponding operator is included in $H^2(\Omega) \times H^2(\Omega)$  and this is one of the conditions for the application of the results in \cite{BK1}.\\

 The {\it isotropic case} when $A(x) = I$ leads to several difficulties. First, the boundary problem (\ref{eq:1.1}) is not parameter-elliptic and the matrix operator related to (1.1) is not self-adjoint. Second, the domain $D(\p)$ of $\p$ is not included in $H^2(\Omega) \times H^2(\Omega)$ (see Section 2) and the results in \cite{BK1}, \cite{BK2} cannot be applied directly. In \cite{R}  Robbiano  proved that in the isotropic case we have
$$N(r) \leq C r^{n/2 + 2}, \: r \to +\infty,$$
covering also the case of complex-valued  $m(x)$.\\

In this paper we study for small $\theta > 0$ the upper bounds of the counting function
$$N(\theta, r)= \#\{ \lambda \in \C: \: \lambda \: {\rm is}\:{\rm (ITE)},\: |\lambda| \leq r,\:  |\arg \lambda| \leq \theta\}.$$ 
 Our main result is the following
\begin{thm} 
Let $0 < \theta < \pi/2$. Then for small $\theta$ there exist $0 < \ep(\theta) \to 0$ as $\theta \to 0$ and $C_{\theta} > 0, \:r(\theta) > 0$ so that
\begin{equation} \label{eq:1.4}
N(\theta, r) \leq  3\sqrt{3}( 1 + \ep(\theta)) \frac{\omega_n}{(2\pi)^n}\int_{\Omega} \Bigl( 1 + (1 + m(x))^{n/2}\Bigr)dx r^{n/2} + C_{\theta},
\end{equation}
provided $r \geq r(\theta).$
\end{thm}

Since outside the angle $\{ z \in \C: \: |\arg z| \leq \theta\}$ we have only finite number of eigenvalues (see \cite{V}, \cite{LV4}), from our result we obtain 
\begin{cor} For the counting function $N(r)$ we have the upper bound
\begin{equation} \label{eq:1.5}
N(r) \leq 3 \sqrt{3} (1 + o(1))\frac{\omega_n}{(2\pi)^n}\int_{\Omega} \Bigl( 1 + (1 + m(x))^{n/2}\Bigr)dx r^{n/2}
\end{equation}
with $o(1) \to 0$ as $r \to \infty.$
\end{cor}
 The constant $3 \sqrt{3}(1 + \ep(\theta))$ is not optimal and the factor $3 \sqrt{3}$ comes from the choice of the  cut-off function $f$ in Section 5 and the Weyl inequality for non-compact operators. By a more complicated argument in Section 5 we can slightly improve this constant. On the other hand, for $n = 1$ and $\Omega = \{|x| \leq 1\}$, $m(x) = const$, a precise Weyl formula has been obtained in \cite{HS}. After the submission of this  paper we have been informed that Robbiano \cite{R1} obtained the leading term in  (\ref{eq:1.5}) without the factor $3 \sqrt{3}(1 + o(1)).$ In the special case $m(x) > 0, \: \forall x \in \bar{\Omega},$ a such asymptotics has been proved recently by Faierman \cite{F}.\\

The problem of sharp Weyl asymptotic of the counting function $N(r)$ with a remainder ${\mathcal O}(r^{n/2 - 1/2})$ is open. Recently, Weyl formula with remainder ${\mathcal O}(r^{n/2 - 1/4})$ has been established in \cite{PV}. In this direction the construction of a semi-classical parametrix in Section 3 and our analysis in Section 5 of the counting function of the self-adjoint operator $Q$ defined below could be useful and we hope to discuss in a further work the problem of the remainder in the asymptotic of $N(r)$.\\

We expect that the (ITE) are mainly concentrated in a neighborhood of the positive real axis and in this direction it is natural to make the conjecture that for every fixed $a > 0$ there exists $\eta = \eta(a) > 0$ such that
\begin{equation}
\# \{\lambda \in \C:\: \lambda \:{\rm is}\:(ITE),\: |\lambda| \leq r,\: |\im \lambda| \geq a\} \sim {\mathcal O}(r^{n/2- \eta}),\: r \to + \infty.
\end{equation}

The proof of Theorem 1 is based on the semi-classical techniques developed in \cite{Sj}, \cite{SjZ1}, \cite{SjZ} for the distribution of the scattering resonances. First we write the spectral problem as a problem for a matrix operator ${\mathcal P}$ which is convenient for our argument. Next our goal is to study the eigenvalues of the operator $h^2 {\mathcal P} - \omega_0,\: 0 < h \ll 1$, where $\omega_0 \in \C$ with $\re \omega_0 =1$ and small $0 < \im \omega_0 < 1$ is fixed.
We apply the strategy of \cite{SjZ} and investigate the eigenvalues of the self-adjoint operator
$$Q = (h^2\pc - \omega_0)^* (h^2\pc - \omega_0),\: 0 < h \ll 1$$
with suitably boundary conditions on $\pa \Omega.$  We show that the operator $Q- z$  is elliptic with parameter for $\Im z \neq 0$ and for $\Re z \leq (\im \omega_0 + \delta)^2, \: |\xi| \geq C_0.$ The next step is to examine the trace class operator $f(Q)$  for $f \in C_0^{\infty}( ]-\infty, (\im \omega_0 + 3\delta)^2[; [0, 1])$ which is equal to 1 on $[0, (\im \omega_0 + 2 \delta)^2]$. To calculate the trace of $f(Q)$, we need a representation of the resolvent $(Q - z)^{-1}$ and for this purpose we construct a semi-classical parameterix for $Q- z$ in the region $|\Im z| \geq h^{\epsilon},\: 0 < \epsilon < 1/2$ and a microlocal parametrix in the region $|\xi| \geq C_0 > 0$ with a large constant $C_0.$ We calculate the leading term of tr $f(Q)$ and  this yields an estimate for the number of the eigenvalues $\mu_1,...,\mu_N$ of the operator $|Q|^{1/2}$ in the interval $[\im \omega_0, \im \omega_0 + 2 \delta].$ Next to obtain a bound for the eigenvalues of $h^2{\mathcal P}$ in the disk $\{z \in \C:\: |z - \omega_0| \leq \im \omega_0 + \delta\}$ we apply the Weyl inequality for non-compact operators proved in Appendix a. in \cite{Sj}. This leads to an estimate of the eigenvalues of ${\mathcal P}$ included in the rectangle $h^{-2}S_0,$ where
$$S_0 = \{ (x, y) \in \R^2:\: 1 - \delta_2 \leq x \leq 1 + \delta_2,\: -X \leq y \leq \im \omega_0\}$$
with $\delta_2 = \sqrt{1- \delta^2} (2 \delta \im \omega_0 + \delta^2)^{1/2},\: \tan \theta = \frac{X}{1 - \delta_2}$ and $X \searrow 0$ as $\delta \searrow 0.$
Summing over the union of rectangles, we obtain the result of Theorem 1.\\

{\bf Acknowledgment.} Thanks are due to Johannes Sj\"ostrand for the discussions concerning the application of Weyl inequality in \cite{Sj}  and to Plamen Stefanov and Evgeny Lakshtanov for many useful discussions and comments. We like also to thank the referees for their remarks and suggestions leading to an improvement of the previous version of the paper.

\section{Elliptic boundary problem for $Q - z$}
\renewcommand{\theequation}{\arabic{section}.\arabic{equation}}
\setcounter{equation}{0}

Consider in $\h = L^2(\Omega) \oplus L^2(\Omega)$ the operator
$$\p = \begin{pmatrix} -\Delta & & 0\\ 0 & & -\fm\Delta\end{pmatrix}$$
with domain $D(\p)$ introduced in Section 1. We have the following
\begin{lem} The operator $\p$ with domain $D(\p)$ is closed in $\h.$
\end{lem}
{\bf Proof.} Let $D(\p) \ni (u_j, v_j) \to (u, v) $ in $\h$ and let $\p (u_j, v_j) \to (f, g)$ in $\h$ as $j \to \infty.$ Since $(\Delta u_j, \Delta v_j) \to (\Delta u, \Delta v)$ in the distribution sense, we deduce $\p (u, v) = (f, g).$ Next, $\Delta (u_j - v_j) \to f - (1+m) g \in L^2(\Omega)$ and
$(u_j - v_j)\big\vert_{\pa \Omega} = \pa_{\nu} (u_j - v_j) \big\vert_{\pa \Omega} = 0.$
The Dirichlet and Neumann problems for $\Delta$ in $\Omega$ are related to closed operators in $L^2(\Omega)$ and this implies 
$$(u - v)\big\vert_{\pa \Omega} = \pa_{\nu} (u - v) \big\vert_{\pa \Omega} = 0.$$ \hfill\qed

It is clear that the (ITE) are the eigenvalues $\lambda \neq 0$ of the operator $\p$. On the other hand, the domain $D(\p)$ is not included in $H^2(\Omega) \oplus H^2(\Omega)$. This leads to difficulties if we want to apply the results in \cite{BK1}, \cite{BK2} for the spectral asymptotics of non-self-adjoint operators with domains included in $H^2(\Omega) \oplus H^2(\Omega).$

To establish the discreteness of the spectrum of $\p$ in $\C \setminus \{0\}$, it is more convenient to pass to the
problem studied by Sylvester \cite{Sy} and Robbiano \cite{R}.  For $\lambda \neq 0$ we reduce  (\ref{eq:1.1}) to the problem
\begin{equation} \label{eq:2.1}
\begin{cases} (\Delta +\lambda(1 + m) )U + mV = 0,\: (\Delta + \lambda) V = 0\: {\rm in}\: \Omega,\\
u = 0,\: \pa_{\nu} u = 0 \: {\rm on}\:\pa \Omega, \end{cases}
\end{equation}
where $U = v - u, V = \lambda u.$
Let 
$$B = \begin{pmatrix}  \Delta_{00} & & m \\ 0 & & \Delta_{- -} \end{pmatrix},\: I_m = \begin{pmatrix} (1 + m) & & 0\\ 0 & & 1 \end{pmatrix}.$$
Here $\Delta_{00}$ is the Laplacian with domain $H^2_0(\Omega)$ and $\Delta_{- -}$ is the Laplacian with domain $\{v \in L^2(\Omega): \: \Delta v \in L^2(\Omega)\}.$ As it is mentioned in \cite{Sy}, $B$ is a closed non-self-adjoint operator in $\h$ with domain
$$D(B) = H^2_0(\Omega) \oplus \{v \in L^2(\Omega): \: \Delta v \in L^2(\Omega)\}.$$
The problem (\ref{eq:2.1}) can be written as $(B +\lambda I_m ) (u, v) = (0, 0).$ In \cite{Sy} it was proved that the resolvent 
$(I_m^{-1} B +\lambda_0)^{-1}$ exists for real $\lambda_0 \ll 1$ and this resolvent is upper triangular compact (UTC) operator 
(see Section 3 in \cite{Sy}) for the definition of (UTC)). Following \cite{Sy} (see also \cite{R} for another proof in a more general setup), $(I_m B +\lambda )^{-1}$ has a meromorphic extension in $\C$. Moreover,  in the Laurent expansion of $(I_m^{-1} B + \lambda )^{-1}$ in a neighborhood of a pole the coefficients are finite dimensional operators. A simple calculus for $\lambda \neq 0$ yields
\begin{equation} \label{eq:2.2}
(\p - \lambda)  = -\begin{pmatrix} 1 & 0 \\ 0 & \fm \end{pmatrix} \begin{pmatrix} 0 & \lambda^{-1}\\ 1 & \lambda^{-1} \end{pmatrix} ( B + \lambda I_m) \begin{pmatrix} -1 & 1\\ \lambda & 0 \end{pmatrix},
\end{equation}
\begin{equation} \label{eq:2.3}
(\p - \lambda)^*  = -\begin{pmatrix} -1 & \bar{\lambda}\\ 1 & 0 \end{pmatrix} (B + \lambda I_m)^* \begin{pmatrix} 0 & 1\\ {\bar{\lambda}}^{-1} & {\bar{\lambda}}^{-1} \end{pmatrix}\begin{pmatrix} 1 & 0 \\ 0 & \fm \end{pmatrix} ,
\end{equation}
\begin{equation} \label{eq:2.4}
(\p - \lambda)^{-1}  = -\begin{pmatrix} 0 & \lambda^{-1}\\ 1 & \lambda^{-1} \end{pmatrix} ( B + \lambda I_m)^{-1}   \begin{pmatrix} -1 & 1\\ \lambda & 0 \end{pmatrix} \begin{pmatrix} 1 & 0 \\ 0 & 1 +m \end{pmatrix}.
\end{equation}
Obviously, we have $\begin{pmatrix} 0 & \lambda^{-1}\\ 1 & \lambda^{-1} \end{pmatrix} \begin{pmatrix} -1 & 1\\ \lambda & 0 \end{pmatrix}= I.$
Consequently, $(\p - \lambda)^{-1}$ is meromorphic in $\C \setminus \{0\}$ and $\p$ has a discreet set of eigenvalues in $\C \setminus \{0\}$ with finite algebraic multiplicity. Moreover, the algebraic multiplicity of the eigenvalues $\lambda_0$ of $\p$ coincide with the algebraic multiplicity of the eigenvalue $-\lambda_0$ of $I_m^{-1}B.$\\

\begin{rem} If $(u, v) \in D(\p)$ is an eigenfunction of $\p$ with eigenvalue $\lambda \neq 0$, then the problem $(\ref{eq:1.1})$ has a week solution. It is proved in Theorem $5.1$, \cite{LV4} that this weak solution is always in $H^2(\Omega) \oplus H^2(\Omega).$ This agrees with the definition of eigenfunction in Section $1.$
\end{rem}

Let $\omega_0 \in \C$ be fixed with $\re \omega_0 > 0, \: \im \omega_0 > 0.$ Consider the operator
$h^2 \p - \omega_0,\: 0 < h \ll 1$.
 Throughout our paper we use the notation $D_{x_j} = \frac{1}{i} \pa_{x_j},\: j = 1,...,n.$ It is important to show that
the boundary problem for $h^2 \p - \omega_0$ is elliptic in the semi-classical sense. 
To study the problem close to the boundary, we pass to local coordinates. Given a point $x_0 \in \pa \Omega$, introduce boundary {\bf normal coordinates} $x = \alpha(y', y_n)$ in a neighborhood of $0$ in $\R^n$ by
$$x = \alpha(y) = \alpha(y', 0) + y_n N(y',0),\: y = (y', y_n),$$
where $N(y', 0)$ is the unit normal at $(y', 0)$ and $y'= (y_1,..., y_{n-1}).$
For simplicity we will use below the notation $x = (x', x_n),\: x'= (x_1,...,x_{n-1})$ for the new coordinates.
Locally near $x_0$, the boundary $\pa \Omega$ is defined by $x_n = 0$ and we have $x_n > 0$ in the domain $\Omega.$ The principal symbol of $-h^2\Delta - \omega_0$ in the new coordinates becomes
$$p_0(x, \xi) = \xi_n^2 + s(x, \xi') - \omega_0,$$

\begin{lem}  The boundary problem
\begin{equation} \label{eq:2.3}
\begin{cases} (h^2D_{x_n}^2 + s(x', 0, \xi') - \omega_0)u = 0,\: x_n > 0,\\
(h^2D_{x_n}^2 + s(x', 0, \xi') - (1+m(x', 0))\omega_0)v = 0, \:x_n > 0,\\
u(0) = v(0),\: (h D_{x_n}u)(0) = (h D_{x_n}v)(0),
\end{cases}
\end{equation}
has only a trivial bounded solution $u = v = 0.$
\end{lem}
{\bf Proof.} 
 Let $\xi_n = \tilde{\lambda}_{\pm}(x',  \xi'),\: j = 1,2,$ be the roots of $\xi_n^2 + s(x',0, \xi')- \omega_0 = 0$ with respect to $\xi_n$ with $\pm \im \tilde{\lambda}_\pm(x', \xi') > 0$ and let $\xi_n = \tilde{\mu}_{\pm}(x', \xi')$ be the roots of $ \xi_n^2 + s(x',0, \xi') -(1 + m(x', 0))\omega_0 = 0$ with respect to $\xi_n$ with $\pm \im\tilde{\mu}_\pm(x', \xi') > 0.$ Then the problem
$$(h^2D_{x_n}^2 + s(x', 0, \xi') - \omega_0)u = 0,  (h^2D_{x_n}^2 +s(x', 0, \xi') -( 1+ m(x))\omega_0)v = 0, \: x_n > 0$$
with  boundary conditions $u(0)= v(0),\: hD_{x_n}u(0) = hD_{x_n}v(0)$ has no bounded non-trivial solutions.  In fact, we have the representation
$$u(x_n) = \tilde{C}_1 (x', \xi') e^{i x_n \tilde{\lambda}_{+}/h},\: v(x_n) = \tilde{C}_2(x', \xi') e^{i x_n \tilde{\mu}_{+}/h}$$
and we get $\tilde{C}_1(x', \xi') = \tilde{C}_2(x', \xi'),\: (\tilde{\lambda}_{+} - \tilde{\mu}_{+}) \tilde{C}_1(x', \xi') = 0.$ It is easy to see that $m(x', 0) \neq 0$ implies $\tilde{\lambda}_{+} \neq \tilde{\mu}_{+}$ and we obtain $\tilde{C}_1 = \tilde{C}_2 = 0.$  \hfill\qed

Let $(h^2 \p - \omega_0)^*$ be the adjoint operator to $h^2 \p - \omega_0$. 
The adjoint operator $B^*$ has the representation
$$B^* = \begin{pmatrix} \Delta_{- -}  &  0\\ m & \Delta_{00}\end{pmatrix},$$
hence the domain of $B^*$ is given by
$$D(B^*) = \{ (u, v) \in L^2(\Omega) \oplus H^2_0(\Omega):\: \Delta u \in L^2(\Omega)\}.$$
Therefore we obtain easily from (2.3) that
the operator $(h^2 \p - \om)^*$ has domain
$$D((h^2\p - \om)^*) = \{ (u, v) \in L^2(\Omega) \oplus L^2(\Omega):\: \Delta u \in L^2(\Omega), \Delta (\fm v) \in L^2(\Omega),\:$$
$$(u + \fm v)\big\vert_{\pa \Omega} = 0,\: \pa_{\nu}(u + \fm v)\big\vert_{\Omega} = 0\}.$$

Next we introduce the operator
$$Q = (h^2\pc - \omega_0)^* (h^2\pc - \omega_0),\: 0 < h \ll 1$$
with boundary conditions
\begin{eqnarray}
\gamma_0 (u - v),\: \gamma_0 \Bigl(\pa_{\nu}( u -v)\Bigr) = 0, \label{eq:2.6} \\
 \gamma_0 (h^2 \Delta u+\omega_0 u) = \gamma_0 \Bigl(\fm\Bigl(-\frac{h^2}{1 + m}  \Delta v - \omega_0 v\Bigr)\Bigr), \label{eq:2.7}\\
\gamma_0 \Bigl( \pa_{\nu}(h^2 \Delta u + \omega_0 u)\Bigr) = \gamma_0 \Bigl(\pa_{\nu} \Bigl(\fm[ - \frac{h^2}{1 + m}\Delta v - \omega_0 v]\Bigr)\Bigr). \label{eq:2.8}
\end{eqnarray}
Clearly, the boundary conditions are determined from the inclusion $$\:(h^2\p - \om)(u, v) \subset D(h^2 \p - \om)^*).$$
 This implies immediately the following

\begin{prop} The operator $Q$ with domain 
$$D = \{(u, v) \in L^2(\Omega) \oplus L^2(\Omega):\: (u, v) \in D(\p), \: (h^2\p - \om)(u, v) \in D( h^2\p - \om)^*)\}$$
 is self-adjoint.
\end{prop}

The determinant of the principal symbol of $Q$ is 
$$q(x, \xi) = (|\xi|^2 - \bar{\omega_0}) (|\xi|^2 - \omega_0)(\fm |\xi|^2 - \bar{\omega}_0) (\fm |\xi|^2 - \omega_0)$$
and the symbol $q(x, \xi)$ is clearly elliptic since
$$|q(x, \xi)| \geq (\im \omega_0)^4.$$

Now we pass to the ellipticity of the boundary problems for $Q$ and $Q - z$.
First we study the boundary problem for $Q$ and establish the following 
\begin{prop}  The boundary problem
\begin{equation} 
\begin{cases} (h^2D_{x_n}^2 + s(x', 0, \xi') - \bar{\omega}_0)(h^2D_{x_n}^2 + s(x', 0, \xi') - \omega_0)u = 0,\: x_n > 0,\\
(h^2D_{x_n}^2 + s(x', 0, \xi') - (1 + m(x', 0))\bar{\omega}_0)(h^2D_{x_n}^2 + s(x', 0, \xi') - (1+m(x', 0))\omega_0)v = 0, \:x_n > 0,\\
u(0) = v(0),\: (h D_{x_n}u)(0) = (h D_{x_n}v)(0),\\
-\Bigl((h^2D_{x_n}^2 + s(x', 0, \xi') - \omega_0)u\Bigr)(0) = \frac{1}{(1 + m(x', 0))^2}\Bigl(\Bigl[h^2D_{x_n}^2 + s(x', 0, \xi') -(1 + m(x', 0)) \omega_0\Bigr] v\Bigr)(0),\\
- hD_{x_n}\Bigl((h^2D_{x_n}^2 + s(x', 0, \xi') - \omega_0)u\Bigr)(0)= \frac{1}{(1 + m(x', 0))^2} hD_{x_n} \Bigl(\Bigl[(h^2D_{x_n}^2 + s(x', 0, \xi')) -(1 + m(x', 0)) \omega_0 \Bigr]v\Bigr)(0).
\end{cases}
\end{equation}
has only a trivial bounded solution $u = v = 0.$
\end{prop}

{\bf Proof.} Multiply the second equation by $\frac{1}{(1 + m(x', 0))^2}$. 
Then 
$$w_1 = - (h^2D_{x_n}^2 + s(x', 0, \xi') - \omega_0)u$$
and 
$$w_2 = \frac{1}{(1 + m(x', 0))^2}\Bigl[h^2D_{x_n}^2 + s(x',0, \xi') -(1 + m(x', 0)) \omega_0\Bigr] v$$
 are solutions of the equations
$$(h^2D_{x_n}^2 + s(x', 0, \xi') - \bar{\omega}_0)w_1 = 0, \: x_n > 0,$$
$$(h^2D_{x_n}^2 + s(x', 0, \xi') - (1+ m(x', 0))\bar{\omega}_0)w_2 = 0, \: x_n > 0.$$
Exploiting the boundary conditions for $w_1$ and $w_2$ and the above argument 
  with $\omega_0$ replaced by $\bar{\omega}_0$, we conclude that $w_1 = w_2 = 0.$ Therefore
$$(h^2D_{x_n}^2 + s(x', 0, \xi') - \omega_0)u = (h^2D_{x_n}^2 + s(x', 0, \xi') - (1 + m(x', 0)) \omega_0) v = 0$$
 Repeating the argument of Lemma 2, we get $u = v = 0$ and this completes the proof. \hfill\qed\\

    Next we choose $ z \in \C \setminus [(\im \omega_0)^2, \infty[$ and consider the operator $Q- z$. We will examine when
$Q-z$ with boundary conditions (\ref{eq:2.6})-(\ref{eq:2.8}) is elliptic in the semi-classical sense. Set $p(x, \xi', \xi_n) = \xi_n^2+ s(x, \xi')$ and let
\begin{equation} \label{eq:2.10*}
(p(x, \xi', \xi_n) - \bar{\omega}_0)(p(x, \xi', \xi_n) - \omega_0) - z = 0
\end{equation}
have roots $\lambda_{j, +}(x, \xi', z),\: j = 1,2$ in the upper half plane and roots $\lambda_{j, -}(x, \xi', z), \: j = 1,2$ in the lower half plane. Also let
\begin{equation} \label{eq:2.11*}
(p(x, \xi', \xi_n) - (1 + m(x))\bar{\omega}_0)(p(x, \xi', \xi_n) - (1+ m(x))\omega_0) - (1 + m(x))^2 z = 0
\end{equation}
have roots $\mu_{j, +}(x, \xi', z),\: j= 1, 2$ in the upper half plane and roots $\mu_{j, -}(x, \xi', z), \: j = 1,2$ in the lower half plane.\\

\begin{lem} If $z \in \C \setminus [(\im \omega_0 + \delta)^2, \infty[$ and $\delta > 0$ is small enough, we have $\lambda_{j, +}(x', 0, \xi', z) \neq \mu_{k, +}(x', 0, \xi', z),\: j = 1,2, \:k = 1,2.$
\end{lem}
 
{\bf Proof.} For simplicity we write $m$ instead of $m(x', 0).$  If we have $\lambda_{j, +} = \mu_{k, +},$ then
$$(p(x', 0, \xi', \lambda_{j, +}) - (1 + m)\bar{\omega}_0)(p(x', 0, \xi', \lambda_{j, +}) - (1+ m)\omega_0) -(1 + m)^2  z = 0$$
yields
$$ - \omega_0 (p(x', 0, \xi', \lambda_{j, +}) - \bar{\omega}_0) - \bar{\omega}_0(p(x', 0, \xi', \lambda_{j, +}) - (1 + m)\omega_0) = (2 + m)z.$$
Thus we get
$$\frac{(m + 2) (|\omega_0|^2 - z)}{2 \re \omega_0} = p(x', 0, \xi', \lambda_{j, +}).$$
On the other hand, the equality $p^2 - 2 \re \omega_0 p + |\omega_0|^2 - z = 0$ with $\xi_n = \lambda_{j, +}$ leads to
$$\frac{( 2 + m)^2(|\omega_0|^2 -z)^2}{4 (\re  \omega_0)^2} - (1 + m)(|\omega_0|^2 - z)  = 0.$$
For small $\delta > 0$ it is clear that $|\omega_0|^2 - z \neq 0.$ In fact, if $z$ is real, we have $(\im \omega_0)^2 - \re z > -o(\delta).$  Next, dividing by $|\omega_0|^2 -z$, we obtain
$$(\im \omega_0)^2 + (\re \omega_0)^2 - z = \frac{4(1 + m)}{(2 + m)^2} (\re \omega_0)^2,$$
hence
$$ (\im \omega_0)^2 - z = -\frac{m^2}{(2 + m)^2} (\re \omega_0)^2$$
which is impossible  for small $\delta.$ \hfill\qed\\

Now we will study the following boundary problem for $Q - z$.
\begin{equation}
\begin{cases}\Bigl[(h^2D_{x_n}^2 + s(x', 0, \xi') - \bar{\omega}_0)(h^2D_{x_n}^2 + s(x', 0, \xi') - \omega_0) - z\Bigr]u = 0,\: x_n > 0,\\
\Bigl[(h^2D_{x_n}^2 + s(x', 0, \xi') - (1+ m(x', 0))\bar{\omega}_0)(h^2D_{x_n}^2 + s(x', 0, \xi') -(1 + m(x', 0)) \omega_0)\\
 -(1 + m(x', 0))^2 z\Bigr] v = 0,\: x_n > 0,\\
 u (0)=  v(0),\: (hD_{x_n}u)(0) = (hD_{x_n}v)(0),\\
-\Bigl((h^2D_{x_n}^2 + s(x', 0, \xi') - \omega_0)u\Bigr)(0) = \frac{1}{(1 + m(x', 0))^2}\Bigl(\Bigl[h^2D_{x_n}^2 + s(x', 0, \xi') -(1 + m(x', 0)) \omega_0\Bigr] v\Bigr)(0),\\
- hD_{x_n}\Bigl((h^2D_{x_n}^2 + s(x', 0, \xi') - \omega_0)u\Bigr)(0)= \frac{1}{(1 + m(x', 0))^2}h D_{x_n} \Bigl(\Bigl[(h^2D_{x_n}^2 + s(x', 0, \xi')) -(1 + m(x', 0)) \omega_0 \Bigr]v\Bigr)(0). \label{eq:2.12}
\end{cases}
\end{equation}

 For simplicity of the notations we will write $m,\: s$ instead of $m(x', 0),\: s(x', 0, \xi').$
We can simplify the last two boundary conditions in (\ref{eq:2.12}). In fact, since $u(0) = v(0)$, the third boundary condition can be written as follows
$$-(1 + m)^2 (h^2D_{x_n}^2 u)(0) - (2 + 2m + m^2) (s - \omega_0)u(0) + m \omega_0u(0)- (h^2D_{x_n}^2 v)(0)= 0.$$
A similar modification can be done for the fourth boundary condition by using $(D_{x_n}u) (0) = (D_{x_n}v)(0)$. Introduce the semi-classical symbol
$$Y(x', \xi'):  = (2 + 2m(x', 0) + m^2(x', 0)) (\omega_0 - s(x', 0, \xi')) + m(x', 0) \omega_0.$$
Then the last two boundary conditions in (\ref{eq:2.12}) have the form
$$-(1 + m)^2 h^2D_{x_n}^2u(0) -h^2 D_{x_n}^2 v(0) + Y u(0) = 0,$$
$$-(1+ m)^2 h^3D_{x_n}^3 u(0) - h^3D_{x_n}^3 v(0) + Y hD_{x_n}u(0) = 0.$$

\begin{prop} For $\im z \neq 0$ the  problem $(\ref{eq:2.12})$ has only the trivial bounded solution $u = v = 0.$ 
\end{prop}
{\bf Proof.} Every bounded solution of (\ref{eq:2.12}) 
satisfies 
\begin{equation} \label{eq:2.13}
D_{x_n}^k u \to 0,\: D_{x_n}^k v \to 0,\: k = 0, 1, 2,3, \: {\rm for}\: x_n \to + \infty.
\end{equation}

We will study the existence of $L^2(\R^+)$-solutions depending on  $(x', \xi')$ of the homogeneous problem (\ref{eq:2.12}) following the approach in \cite{SjZ}. Let 
$$U = (u(x', x_n, \xi', z), v(x', x_n, \xi', z)) \in H^4 ([0, +\infty[: \C^2),$$
 where we consider $(x'$, $\xi')$ and $z$ as parameters. Assume that $U \in H^4([0, +\infty[: \C^2)$ is a solution of (\ref{eq:2.12}) such that 
(\ref{eq:2.13}) holds. Let
$${\bf p}(x', 0, \xi', hD_{x_n}): = \begin{pmatrix} h^2D_{x_n}^2 + s(x', 0, \xi')- \omega_0 & & 0\\
0 & & \frac{1}{(1 + m(x', 0))^2} \Bigl(h^2D_{x_n}^2 + s(x', 0, \xi') - (1 + m(x', 0)) \omega_0\Bigr) \end{pmatrix}.$$

By using the boundary conditions in (\ref{eq:2.12}), we obtain
\begin{equation} \label{eq:2.14}
( ({\bf p}^*{\bf p}- z) U, U) = - z \|U\|^2 + \|{\bf p} U\|^2,
\end{equation}
where $(. ,.)$ is the scalar product in $L^2(\R^+: \C^2)$ and $\|.\|$ is the $L^2(\R^+: \C^2)$- norm.  Therefore
\begin{equation} \label{eq:2.15}
 \|U\| \leq \frac{1}{|\im z|}\|(Q(x', 0, \xi', hD_{x_n}) - z) U\|.
\end{equation}
This proves that for fixed $(x', \xi')$ we have only a trivial bounded solutions of  (\ref{eq:2.12}).
\hfill\qed\\

For $\Im z \neq 0$, the estimate (\ref{eq:2.15}) holds for every $U \in C_0^{\infty}([0, +\infty[: \C^2)$ satisfying the boundary problem (\ref{eq:2.12}) with $h = 1$ since its proof is based on energy estimates. Thus we deduce that for $U = (u, v) \in H^4(]0, +\infty[:\C^2)$ satisfying the boundary problem in (\ref{eq:2.12}) with $h = 1$ we have
\begin{equation} \label{eq:2.16}
\|U\|_4 \leq \frac{C}{|\im z|} \Bigl\|\Bigl(Q(x', 0, \xi', D_{x_n}) - z\Bigr) U\Bigr\|,
\end{equation}
where the constant $C$ depends of $(x', \xi')$. For $|x'| \leq \rho, |\xi'| \leq C_0$ we may take $C$ uniformly with respect to $(x', \xi').$ Here and below $\|.\|_k$ denotes the norm in the Sobolev space $H^k(]0, +\infty[: \C^2).$\\ 

Applying (\ref{eq:2.16}),  it is easy to obtain an estimate for the  bounded solutions of the non-homogeneous problem
\begin{equation} \label{eq:2.17}
\begin{cases}\Bigl[(D_{x_n}^2 + s(x', 0, \xi') - \bar{\omega}_0)(D_{x_n}^2 + s(x', 0, \xi') - \omega_0) - z\Bigr]u = F_1,\: x_n > 0,\\
\Bigl[(D_{x_n}^2 + s(x', 0, \xi') - (1+ m(x', 0))\bar{\omega}_0)(D_{x_n}^2 + s(x'0, \xi') -(1 + m(x', 0)) \omega_0)\\
 -(1 + m(x', 0))^2 z\Bigr] v = F_2,\: x_n > 0,\\
 u (0) -  v(0) = f_0,\: (D_{x_n}u)(0) - (D_{x_n}v)(0) = f_1,\\
-\Bigl((D_{x_n}^2 + s(x', 0, \xi') - \omega_0)u\Bigr)(0) - \frac{1}{(1 + m(x', 0))^2}\Bigl(\Bigl[D_{x_n}^2 + s(x', 0, \xi') -(1 + m(x', 0)) \omega_0\Bigr] v\Bigr)(0) = f_2,\\
- D_{x_n}\Bigl((D_{x_n}^2 + s(x', 0, \xi') - \omega_0)u\Bigr)(0)\\- \frac{1}{(1 + m(x', 0))^2} D_{x_n} \Bigl(\Bigl[(D_{x_n}^2 + s(x', 0, \xi')) -(1 + m(x', 0)) \omega_0 \Bigr]v\Bigr)(0) = f_3.
\end{cases}
\end{equation}
Let $F = (F_1, F_2).$ Choose functions $H_j(x_n) \in C_0^{\infty}([0, + \infty[),\: j = 0,1,2,3,$ so that
$$ D_{x_n}^k H_j(0) = \delta_{k, j},\: k, j = 0,1,\: D_{x_n}^k H_j (0) = \delta_{k, j}, k = 0,1,2,3,\: j = 2, 3.$$

Consider
$$\tilde{u} = u - f_0 H_0,\: \tilde{v} = v - f_1 H_1.$$
It is clear that
$$\tilde{u}(0) - \tilde{v}(0) = 0,\: D_{x_n} \tilde{u}(0) - D_{x_n}\tilde{v}(0) = 0.$$
Then changing $F, f_2$ and $f_3$, we may write (\ref{eq:2.17}) as a boundary problem for $(\tilde{u}, \tilde{v})$ with two homogeneous boundary conditions. Now let $Y= Y(x', \xi')$ be the symbol introduced above. 
The last two boundary conditions for the problem for $(\tilde{u}, \tilde{v})$ can be written in the form 
$$-(1+ m)^2 D_{x_n}^2 \tilde{u}(0) - D_{x_n}^2 \tilde{v}(0) + Y(x', \xi')\tilde{u}(0) = f_4,$$
$$-(1+m)^2 D_{x_n}^3 \tilde{u}(0) - D_{x_n}^3 \tilde{v}(0) + Y D_{x_n} \tilde{u}(0) = f_5,$$
where $f_4, f_5$ can be expressed by $f_0, f_1, f_2, f_3, H_0, H_1.$ Introduce
$$u_1 = \tilde{u} + f_4(1 + m(x', 0))^{-2}H_2 + f_5(1 + m(x', 0))^{-2} H_3 ,\:v_1 = \tilde{v}.$$
Then
$$-(1 + m(x', 0))^2 D_{x_n}^2 u_1(0) - D_{x_n}^2 v_1(0) + Y u_1(0) = f_4 - f_4 = 0$$
and the last boundary condition for $(u_1, v_1)$ is also homogeneous.
Thus $(u_1, v_1)$ satisfies a boundary problem with homogeneous boundary conditions. On the other hand,
$$\Bigl(Q(x', 0, \xi', D_{x_n}) - z\Bigr) (u_1, v_1) = F + G = \tilde{F},$$
where $G$ depends on $F, f_0, f_1, f_2, f_3, H_0, H_1, H_2, H_3.$ We obtain
$$\|\tilde{F}\| \leq \|F\| + C_1 \sum_{j=0}^3 |f_j|$$ 
and, applying the estimate (\ref{eq:2.15}), we deduce that the bounded solution $U = (u, v)$ of the problem (\ref{eq:2.17}) satisfies
\begin{equation} \label{eq:2.18}
\|U\|_4 \leq \frac{C_2}{|\im z|} \Bigl(\|F\| + \sum_{j=0}^3 |f_j|\Bigr).
\end{equation}

In the following we assume that $z \in \C \setminus [(\im \omega_0 + \delta)^2, + \infty[$ with small $\delta > 0$ and $|z| \leq a_0$.  We have proved that the problem (\ref{eq:2.12}) does not have non trivial  bounded solutions only for $\Im z \neq 0.$ The case when $z \leq 0$  is real or $z$ is small enough can be treated exploiting (\ref{eq:2.14}). The analysis of other cases when $z$ is real is more complicated.  However, we will show below that for $|\xi'| \geq C_0$ with a large fixed constant $C_0 > 0$ the boundary problem (\ref{eq:2.12}) does not have non trivial bounded solutions for $ \re z \leq (\im \om + \delta)^2$. This argument plays a crucial role in the construction of a semi-classical parametrix in the region $(A)$ in Section 3 as well as in the calculation of the trace of $f(Q)$ in Section 4.\\

We start with the case  $ r = \re z - (\im \omega_0)^2 \neq 0$. It is clear that $\lambda_{1, +} \neq \lambda_{2, +}$ and $\mu_{1, +} \neq \mu_{2, +}$. Then a bounded solution of (\ref{eq:2.12}) has the form
\begin{equation} \label{eq:2.18}
\begin{cases}
u(x', x_n, \xi'; z) = C_1(x', \xi') e^{i x_n \lambda_{1, +}/h } + C_2(x', \xi') e^{i x_n \lambda_{2, +}/h},\\
v(x', x_n, \xi'; z) = C_3(x', \xi') e^{ix_n \mu_{1, +}/h} + C_4(x', \xi') e^{i x_n \mu_{2, +}/h}.
\end{cases}
\end{equation}

For the coefficients $C_j(x', \xi'),\: j = 1,...,4,$ in the representation of $u(x', x_n, \xi'; z), v(x', x_n, \xi'; z)$  we obtain a linear system with matrix
\begin{equation} \nonumber
\Lambda(x', \xi', z)  = \begin{pmatrix} 1 & & 1 & & - 1 & & - 1 \\
\lambda_1 & & \lambda_2 & & -\mu_1 & & -\mu_2\\
-(1 + m)^2\lambda_1^2 + Y & & -(1+m)^2\lambda_2^2 + Y& & - \mu_1^2 & & -\mu_2^2\\
-(1+m)^2\lambda_1^3  + \lambda_1 Y& & -(1+ m)^2\lambda_2^3  + \lambda_2 Y& &  -\mu_1^3 & & -\mu_2^3
 \end{pmatrix},
\end{equation}
where for simplicity we write $\lambda_{j},\: \mu_k$ instead of $\lambda_{j, +}, \: \mu_{k, +}$ .
Clearly,
\begin{equation} \nonumber
\frac{\det\Lambda}{ (\lambda_2 - \lambda_1)} = \det \begin{pmatrix} 1 & & \lambda_2 - \mu_1 & & \mu_1 - \mu_2\\
-(1+m)^2(\lambda_1 + \lambda_2) & & -(1+m)^2 \lambda_2^2  - \mu_1^2 +Y& &  \mu_1^2- \mu_2^2  \\
-(1+m)^2(\lambda_1^2 + \lambda_1 \lambda_2 + \lambda_2^2)+Y & & -(1+ m)^2\lambda_2^3 - \mu_1^3 +\lambda_2 Y & & \mu_1^3 - \mu_2^3\\
\
\end{pmatrix}. 
\end{equation}
Multiplying the first column by $\lambda_2$ and taking the difference with the second one, we get
$$\det \Lambda (x', \xi', z) = -(\lambda_2 - \lambda_1)  (\mu_2 - \mu_1) \det \Lambda_1(x', \xi', z),$$
where
\begin{eqnarray*}
 \Lambda_1 : =
 \begin{pmatrix} 1  & & -\mu_1 & & 1\\
-(1+m)^2 (\lambda_1 + \lambda_2) & & (1+m)^2 \lambda_1\lambda_2 - \mu_1^2+ Y & & \mu_1 + \mu_2 \\
-(1+m)^2(\lambda_1^2 + \lambda_1 \lambda_2 + \lambda_2^2) + Y & & (1+ m)^2\lambda_1\lambda_2(\lambda_1 + \lambda_2) - \mu_1^3 & &  \mu_1^2 + \mu_1\mu_2 + \mu_2^2 \end{pmatrix}.
\end{eqnarray*}
Next multiplying the third column by $\mu_1$ and adding it to the second one, we get
$$\det \Lambda_1(x', \xi', z) $$
$$= \det \begin{pmatrix} 1  & & 0 & & 1\\
-(1+m)^2 (\lambda_1 + \lambda_2) & & (1+m)^2 \lambda_1\lambda_2 + \mu_1\mu_2+ Y & & \mu_1 + \mu_2 \\
-(1+m)^2(\lambda_1^2 + \lambda_1 \lambda_2 + \lambda_2^2) + Y & & (1+ m)^2\lambda_1\lambda_2(\lambda_1 + \lambda_2) + \mu_1\mu_2(\mu_1 + \mu_2) & &  \mu_1^2 + \mu_1\mu_2 + \mu_2^2 \end{pmatrix}
$$ 

$$= (1 + m)^2 \Bigl( \mu_1\mu_2 (\lambda_1 + \lambda_2)^2 - (\lambda_1\lambda_2+ \mu_1 \mu_2) (\mu_1 + \mu_2)(\lambda_1 + \lambda_2) + \lambda_1 \lambda_2(\mu_1 + \mu_2)^2\Bigr)$$
$$-\Bigl( (1+m)^2 \lambda_1 \lambda_2 + \mu_1\mu_2\Bigr)^2 + \Bigl((1+m)^2(\lambda_1^2 + \lambda_2^2) + (\mu_1^2 + \mu_2^2)\Bigr)  Y - Y^2. $$
We will prove the following
\begin{prop} Assume
\begin{equation} \label{eq:cond}
(2 + 3m + m^2) \im \omega_0 \leq \eta m(m+1) \re \omega_0, \: 0 < \eta < 1.
\end{equation}
 Then for  $ z \in \C \setminus \Bigl([(\im \omega_0 + \delta)^2, + \infty[ \cup \{(\im \omega_0)^2\}\Bigr)$, $|z| \leq a_0$ and $\delta > 0$ small enough there exists a constant $C_0 > 0$ such that for $|\xi'| \geq C_0$ we have the estimate
\begin{equation} \label{eq:2.19}
|\det \Lambda_1(x', \xi', z)| \geq C_2 
\end{equation}
with $C_2 > 0$ independent of $(x', \xi')$ and $z$. In particular, for these values of $z$ and $\xi'$ the problem $(\ref{eq:2.12})$ has only a trivial bounded solution $u = v = 0.$
\end{prop}

{\bf Proof.} 
We start with  a representation of $\lambda_{j, +}$ and $\mu_{k, +}$ for $|\xi'| \geq C_0$. If we have
$$p(x', 0, \xi', \xi_n)^2 - 2 (\re \omega_0) p(x', 0, \xi', \xi_n) + |\omega_0|^2 - z = 0,$$ 
then $p(x', 0, \xi', \xi_n) = \re \omega_0 \pm \sqrt{ \re z - (\im \omega_0)^2 + i \im z}.$ Set 
$$r = \re z - (\im \omega_0)^2 \in \R,\: \im z = w.$$

We choose the branch $0 \leq {\rm Arg} \: z < 2\pi.$  Then for $|\xi'|\geq C_0$ and large $C_0 > 0$ we have
$$\lambda_{1, \pm} = \pm\sqrt{\re \omega_0 - s + \sqrt{r + iw}},\: \lambda_{2, \pm} = \pm\sqrt{\re \omega_0 - s - \sqrt{r + i w}},$$
$$\mu_{1, \pm} = \pm\sqrt{(m+1) \re \omega_0 - s + (1 + m)\sqrt{r + i w}},$$
$$ \mu_{2, \pm} = \pm\sqrt{(m+1) \re \omega_0 -s -(1 + m)\sqrt{r + i w}}.$$

It is clear that the equations (\ref{eq:2.10*}), (\ref{eq:2.11*}) imply
\begin{equation*} 
\lambda_1^2 + \lambda_2^2 = 2(\re \omega_0-s),\:\mu_1^2 + \mu_2^2 = 2(1 + m) \re \omega_0 - 2s.
\end{equation*}
Therefore,
$$(1+m)^2 (\lambda_1^2 + \lambda_2^2) + (\mu_1^2 + \mu_2^2)$$ 
$$= 2\Bigl[ ((1+m)^2 + (1+m)) \re \omega_0 - ((1+m)^2 +1)s\Bigr] = 2 \re Y$$
and we obtain\
$$[ 2 \re Y - Y] Y = (\re Y - i \im Y) Y = |Y|^2.$$
Set $\zeta = \sqrt{r + iw}$. Then for $|\xi'| \to + \infty$, we have 
$$\lambda_{1,2} = i\sqrt{s}  - i\frac{ \re \omega_0  \pm\zeta}{2\sqrt{s}} - i\frac{(\Re \omega_0 \pm\zeta)^2}{8 s^{3/2}} + \Oc(\frac{1}{s^{5/2}}),$$
$$\mu_{1, 2} = i\sqrt{s}  - i\frac{ (m+1)(\re \omega_0)  \pm\zeta}{2\sqrt{s}} - i\frac{((1+ m)\Re \omega_0  \pm \zeta)^2}{8 s^{3/2}} + \Oc(\frac{1}{s^{5/2}}),$$
where we have the sign $(+)$ for $\lambda_1, \mu_1$ and the sign $(-)$ for $\lambda_2, \mu_2.$ This implies
$$\lambda_1 \lambda_2 = -s + \re \omega_0 +\frac{r + i w}{2 s}+ \Oc(\frac{1}{s^2})$$
and 
$$\mu_1 \mu_2 =  -s + (m + 1)\re \omega_0 + \frac{(m+1)^2 (r + i w)}{2s} +{\mathcal O}(\frac{1}{s^2}).$$
Also notice  that we have
$$\lambda_1 + \lambda_2 = 2i \sqrt{s} - \frac{i \re \omega_0}{\sqrt{s}} + {\mathcal O}(\frac{1}{s^{3/2}}),$$
$$\mu_1 + \mu_2 = 2 i \sqrt{s} -\frac{i (m + 1) \re \omega_0}{\sqrt{s}} + {\mathcal O}(\frac{1}{s^{3/2}}).$$

Therefore
$$(1+m)^2 \lambda_1 \lambda_2 + \mu_1\mu_2 
= -(2 + 2m + m^2) s + (2 + 3m + m^2) \re \omega_0 +  \frac{(m + 1)^2(r + i w)}{s}+ \Oc(\frac{1}{s^2})$$
$$= \re Y + \frac{(m + 1)^2(r + i w)}{s}+ \Oc(\frac{1}{s^2})$$
and we deduce
$$\Bigl((1+m)^2 \lambda_1 \lambda_2 + \mu_1\mu_2\Bigr)^2 =  |\re Y|^2 - 2(2 + 2m + m^2)(m+ 1)^2(r + i w) + \Oc(\frac{1}{s}).$$
We multiply by (-1) the above expression, sum it with $|Y|^2$ and obtain
\begin{eqnarray}
 -\Bigl((1+m)^2 \lambda_1 \lambda_2 + \mu_1\mu_2\Bigr)^2 + |Y|^2 = |\im Y|^2 + 2(2 + 2m + m^2)(m+ 1)^2(r + i w) + \Oc(\frac{1}{s})\nonumber \\
= (2 + 3 m + m^2)^2 (\im \omega_0)^2 +  2(2 + 2m + m^2)(m+ 1)^2(r + i w) + \Oc(\frac{1}{s}). \label{eq:2.20}
\end{eqnarray}

It remains to examine the behavior of the term
$$Z : = \mu_1\mu_2 (\lambda_1 + \lambda_2)^2 - (\lambda_1\lambda_2 + \mu_1 \mu_2) (\mu_1 + \mu_2)(\lambda_1 + \lambda_2) + \lambda_1 \lambda_2(\mu_1 + \mu_2)^2$$
$$= \Bigl((\lambda_1 + \lambda_2) - (\mu_1 + \mu_2)\Bigr) \Bigl( \mu_1 \mu_2(\lambda_1 + \lambda_2) - \lambda_1 \lambda_2(\mu_1 + \mu_2)\Bigr).$$
We have
$$(\lambda_1 + \lambda_2) - (\mu_1 + \mu_2) = \frac{i m \re \omega_0}{\sqrt{s}} + \Oc(\frac{1}{s^{3/2}}),$$

$$\mu_1 \mu_2(\lambda_1 + \lambda_2) - \lambda_1 \lambda_2(\mu_1 + \mu_2) = i m \re \omega_0 \sqrt{s} +  \Oc(\frac{1}{\sqrt{s}}),$$
hence
\begin{equation} \label{eq:2.21}
(1+m)^2Z = - m^2 (1+m)^2(\re \omega_0)^2 + \Oc(\frac{1}{s}).
\end{equation}
Taking account of (\ref{eq:2.20}), (\ref{eq:2.21}), for $|\im z| \geq \eta_0 > 0$ we obtain an uniform lower bound for $|\det \Lambda_1|$ since
$\im (\det \Lambda_1) = 2( 2 + 2m + m^2)(m+1)^2 \im z + \Oc(\frac{1}{s}).$ To treat the case $|\im z| < \eta_0$, we exploit the condition (\ref{eq:cond}) which yields
$$(2 + 3m + m^2)^2 (\im \omega_0)^2 - \eta^2 \:m^2(m+1)^2 (\re \omega_0)^2 \leq 0.$$ 

Since $r \leq C\delta,$ the leading term of the real part of $\det \Lambda_1$ is less than
$$- (1- \eta^2)m^2(m+1)^2 (\re \omega_0)^2 +2(2 + 2m + m^2) (m+1)^2 C \delta.$$
For small $\delta > 0$ depending on $\re \omega_0$ and $m$ we obtain an uniform lower bound
 and this completes the proof. \hfill\qed\\

\begin{rem} For $\im z \neq 0$ the lower bound of $|\det \Lambda_1|$ has the form $C(m) |\im z|$ and this agrees with the estimates
$(2.15), (2.18)$. For our purpose Proposition $4$ is sufficient since in Section $5$ we will choose $\re \omega_0 = 1$ and $0 < \im \omega_0 < 1$ small.
\end{rem}

Obviously, for $z = (\im \om)^2$ and $|\xi'| \geq C_0$ we have double roots 
$$ \lambda_{+}: = \sqrt{\re \om - s} ,\:\mu_{+} = : \sqrt{(1 + m)\re \om -s}.$$
In this case a bounded solution of (\ref{eq:2.12}) has the form
\begin{equation*} 
\begin{cases}
u(x', x_n, \xi'; z) = C_1(x',\xi') e^{i x_n \lambda_{+}/h} + C_2(x', \xi') \frac{i x_n}{h} e^{i x_n \lambda_{+}/h},\\
v(x', x_n, \xi'; z) =  C_3(x',\xi') e^{i x_n \mu_{+}/h} + C_4(x', \xi') \frac{i x_n}{h} e^{i x_n \mu_{+}/h}.
\end{cases}
\end{equation*}
For the coefficients $C_j(x', \xi')$ we obtain a linear system with the matrix
\begin{equation} \nonumber
\Lambda_2(x', \xi', z): = \begin{pmatrix} 1 & & 0 & & - 1 & & 0 \\
\lambda & & 1 & & -\mu & & -1 \\
-(1 + m)^2 \lambda^2 + Y & & -2(1 + m)^2\lambda & & -\mu^2 & & - 2\mu \\
-(1 +m)^2 \lambda^3 + Y \lambda & & -3(1 + m)^2\lambda^2 + Y& & - \mu^3 & & -3\mu^2 \end{pmatrix},
\end{equation}
where we write simply $\lambda, \mu$ for $\lambda_{+}, \: \mu_{+}$.  A calculation of $\det \Lambda_2$ implies
$$\det \Lambda_2 = \mu^4 +(1+m)^4 \lambda^4 + (1 + m)^2( 4 \lambda \mu^3 + 4 \lambda^3 \mu- 6 \lambda^2 \mu^2)- 2((1 +m)^2\lambda^2 + \mu^2)Y + Y^2$$
$$= ((1+m)^2 \lambda^2 + \mu^2)^2 + 4(1+m)^2 \lambda \mu( \lambda - \mu)^2 - 2((1+m)^2\lambda^2 + \mu^2) Y + Y^2.$$
$$= [(1+m)^2\lambda^2 + \mu^2 - Y]^2 + 4(1+m)^2 \lambda \mu( \lambda - \mu)^2.$$
We have 
$$(1+m)^2 \lambda^2 + \mu^2 -Y$$
$$= (1+m)^2 (\re \omega_0 -s) + (1+m) \Re \omega_0 -s + (2 + 2m + m^2) s -(2 + 3 m + m^2) \omega_0$$
$$= -(2 + 3m + m^2) i\im \omega_0.$$
The term $4\lambda \mu (\lambda- \mu)^2$ is the analog of the term $-Z$ in the proof of Proposition 4. By using the
continuity of the roots $\lambda_{j, +}, \: \mu_{j, +}$ as $z \to (\im \om)^2,$ we get
$$4 (m+1)^2\lambda \mu(\lambda- \mu)^2 = m^2 (m+1)^2(\re \omega_0)^2 + \Oc(\frac{1}{s}).$$
Thus we obtain the following
\begin{prop}  Assume the condition $(\ref{eq:cond})$ fulfilled. Then for $z = (\im \omega_0)^2$ there exists a constant $C_0 > 0$ such that for $|\xi'| \geq C_0$ we have the estimate
\begin{equation} \label{eq:2.22}
|\det \Lambda_2 (x', \xi', z)| \geq C_3 
\end{equation} 
with $C_3 > 0$ independent of $(x', \xi')$. In particular, for these values of $\xi'$ the problem $(2.12)$ has only a trivial bounded solution $u = v = 0.$
\end{prop}
\begin{rem} We can obtain $\Lambda_2$ from $\frac{\Lambda}{(\lambda_2 - \lambda_1)(\mu_2- \mu_1)} = -\Lambda_1$ by taking the limit $z \to (\Im \omega_0)^2.$  This explains the appearance of the sign $(-)$ in the expression of $\Lambda_2.$
\end{rem}

\section{Semi-classical parametrix}

In this section we assume that $z$ belongs to a compact set $Z_0$ in $\C \setminus [((\im \omega_0) + \delta)^2, \infty[$ which is independent on $h$. We assume also that $\delta > 0$ is small enough and $\omega_0$ is such that (\ref{eq:cond}) holds.
This guarantees that the estimates (\ref{eq:2.19}) and (\ref{eq:2.22}) are satisfied for $|\xi'| \geq C_0.$ Throughout the paper we use for h-pseudodifferentrial operator the standard quantization (see \cite{DS}, \cite{SjZ}).

Since the boundary problem (\ref{eq:2.12}) is not elliptic for all $z \in Z_0$ and all $\xi'\in \R^{n-1}$, we will construct a semi-classical parametrix in two regions:
$$(A):\: z\in Z_0, \: |\xi'| \geq C_0,$$
$$(B):\: z \in Z_0,\: |\im z | \geq h^{\epsilon}.$$ 
Here $C_0 > 0$ is the constant from Propositions 4, 5 and $0 \leq \epsilon < 1/2$ is a small fixed constant. In the region (A) we construct a microlocal parametrix which will be defined precisely below and the condition $|\xi'| \geq C_0$ must be interpreted in local coordinates
$(x', \xi').$
 The semi-classical parametrix will be a sum of an operator with matrix symbol $R(x, \xi; z, h)$ defined for all $(x, \xi) \in \R^{2n}$ and a matrix-valued operator corresponding to the boundary conditions. \\

 We extend $m(x)$ to $\R^n$ so that 
$m(x) \in C_0^{\infty}(\R^n)$ with $m(x) = 0$ for $|x| \geq \rho > 0.$ The operator $Q$ with coefficient $m(x)$ extended for all $x \in \R$ will be denoted also by $Q$. First we consider the region (B) and the modifications for (A) will be explained later. We start with the construction of a parametrix for $Q - z$ in $\R^n$. We assume $|\im z| \geq h^{\epsilon}$ and we search a matrix symbol $R(x, \xi; z, h)$ determined for all $x \in \R^n$ and all $\xi \in \R^n$ so that 
$$ (Q - z) \#_h R \sim I.$$
Here $a\#_h b$ denotes the composition of semi-classical symbols $a$ and $b$ (see Section 7 in \cite{DS} and Appendix A in \cite{HK3}) and we use the same notation for the operator $Q$ and for its symbol. We denote by $Op_h(c)$ the h-pseudo-differential operator with semi-classical symbol $c$.
Set
$$Q_1(x, hD_x) =(-h^2\Delta  - \bar{\omega}_0)(-h^2 \Delta - \omega_0),$$
$$ Q_2(x, hD_x)=(-h^2\Delta  \fm - \bar{\omega}_0)(-h^2 \fm \Delta - \omega_0).$$
Denote by $q_k(x, \xi)$ the principal symbols of the semi-classical operators $Q_k(x, hD_x), k = 1,2.$ Clearly,  the symbol $q_1$ is independent on $x$ but in the following we use the notation $q_1(x, \xi)$ which is convenient when we treat the boundary conditions. 
Since the operator $Q- z$ has a diagonal form, the construction follows a standard argument and we find $R$ as a matrix-valued operator
$$R = \begin{pmatrix} R_1 & & 0\\
0 & & R_2 \end{pmatrix},$$
so that 
$$(Q_1- z) R_1 = I + K_1,\:(Q_2- z) R_2 = I + K_2.$$
Moreover,  for the distribution kernels $K_j(x, y; z, h)$ of the operators $K_j, j = 1,2,$ we have
$$|\partial_x^{\alpha} \partial_y^{\beta} K_j(x, y; z, h)| \leq C_{\alpha, \beta, N} h^N,\:\: \forall \alpha, \forall \beta, \forall N.$$

Below we use the spaces of symbols $S_{\epsilon}^{m, k}:=S_{\epsilon}^{m, k}(\R^{2n})$ (see \cite{SjZ} and Appendix A in \cite{HK3} for the properties of the operators with symbols in these classes.) For reader convenience recall that  $a(x, \xi; h) \in S_{\epsilon}^{m, k}(\R^{2n})$ if $a$ is $C^{\infty}$ smooth with respect to $(x, \xi) \in \R^n \times \R^n$ and such that\\
(i) for $\xi$ outside some $h-$ independent compact set $a$ satisfies $a \in S^k(\R^{2n})$, where $S^k(\R^{2n})$ denotes the H\"ormander space, i.e.,  $\vert \partial_x^\alpha \partial_\xi^\beta a(x,\xi; h)\vert \leq C_{\alpha,\beta}(1+\vert \xi\vert)^{k-\vert \beta\vert  }$ uniformly with respect to $h$,\\
(ii) For every compact set $\Xi$ in $\R^n$ we have the estimates
$$|\pa_x^{\alpha} \pa_{\xi}^{\beta} a(x, \xi ;h )| \leq C_{\alpha, \beta, \Xi} h^{-m} h^{-\epsilon(|\alpha| + |\beta|)},\: \forall (x, \xi) \in \R^n \times \Xi,\: \forall \alpha, \forall \beta \in \N^n.$$
For a matrix-valued symbol $A$ we say that $A \in S_{\epsilon}^{m, k}(\R^n; {\mathcal M}_2(\C))$ if every component $a_{i,j}$ of the matrix $A$ is in $S_{\epsilon}^{m, k}.$\\

It is clear that
$$\big|\Bigr(|\xi|^2 - \bar{\omega_0}) (|\xi|^2 - \omega_0)- z\Bigr)\Bigr((\fm |\xi|^2 - \bar{\omega}_0) (\fm |\xi|^2 - \omega_0)- z\Bigr)\big| \geq |\im z|^2. $$
Using the construction of a semi-classical parametrix  depending on a parameter in the region (B) (see for example, \cite{DS}, \cite{SjZ}, \cite{HK3}),  we  find ${\mathcal R}_N = \sum_{j=0}^N h^jc_j$ with
$c_j \in S_{\epsilon}^{\epsilon (1 +j), - 4 - j}(\R^{2n}; {\mathcal M}_2({\mathbb C}))$ so that
$$(Q- z)Op_h({\mathcal R}_N) = I +  h^{N+1}Op_h (c_{N+1}).$$

The symbol of the operator $R$ is defined as $R \sim \sum_{j=0}^{\infty} h^j c_j. $ For example, the symbol $R_1$ has
an asymptotic expansion
$$R_1 \sim \frac{1}{q_1 - z} +h \frac{b_1(x, \xi, z)}{(q_1 - z)^3} + h^2\frac{b_2(x, \xi, z)}{(q_1 - z)^5} + ...$$
where $b_j(x, \xi, z)$  is a polynomial in $z$ with smooth coefficients (see Section 8 in \cite{DS} for more details).

Since the construction is a repetition of that in  \cite{DS}, \cite{SjZ}, \cite{HK3} we omit the details.\\

Now we pass to the construction close to the boundary. Our argument is similar to that in Section 3 in \cite{SjZ} and Section 4 in \cite{HK3}. We introduce local {\bf geodesic coordinates}
 in a neighborhood of a point $x_0 \in \partial \Omega$ as in the previous section. For simplicity we will use for variables the notation $(x', x_n, \xi', \xi_n)$.

Let $\lambda_{j, \pm}(x, \xi', z)$ and $\mu_{j, \pm}(x, \xi', z),\: j = 1,2,$ be the roots defined in the previous section.   We can choose  $R \in S_{\epsilon}^{\epsilon, - 4}$ with a holomorphic extension in $\xi_n$-variable in the domain
$$\Omega_1(x, \xi', z) = \{\xi_n \in \C,\:|\xi_n| \leq \eta^{-1} \la \xi'\ra,\: |\xi_n - \lambda_{j, \pm} | > \eta \la \xi'\ra,\: |\xi_n - \mu_{j, \pm}| > \eta\la \xi'\ra, \: j = 1,2\}$$ 
for arbitrary fixed $\eta > 0.$ Here and below we use the notation $\la \xi'\ra = (1 + |\xi'|^2)^{1/2}$ and we will denote by $C$ some positive constants which may change from line to line.\\

 For $|\xi'| \leq C_0$ we have
$$|\im \lambda_{j, \pm}(x, \xi', z) | \geq h^{\epsilon}/C,\: |\im \mu_{j, \pm}(x, \xi', z)| \geq h^{\epsilon}/C, j = 1,2.$$
Consider the domain
$$\Omega(x, \xi', z) = (\Omega_1 \cap \{\xi_n \geq \eta \la \xi'\ra, \:|\xi'| \geq C_0\}) \cup \Omega_2 \cup \Omega_3,$$ 
where
$$\Omega_2 (x, \xi', z)= \{|\xi'| \leq C_0,\: |\xi_n| \leq C, 0 \leq \im \xi_n \leq \frac{h^{\epsilon}}{C}\},$$
$$\Omega_3(x, \xi', z) = \{|\xi'| \leq C_0, \:|\xi_n| \leq C, \: \im \xi_n \geq \frac{h^{\epsilon}}{C},\: |\xi_n - \lambda_{j, +}| \geq C^{-1},\: |\xi_n - \mu_{j, +}| \geq C^{-1},\: j = 1,2\}.$$
We have a holomorphic extension of $R$ with respect to $\xi_n$ in $\Omega(x', \xi', z).$ Let 
$$\gamma = \gamma(x', \xi') \subset \Omega (x, \xi', z)$$
 be a simple curve which encircles $\lambda_{j, +}$ and $\mu_{j, +}$, $j = 1, 2,$ in  a positive sense.
Introduce the operators $\Pi_j: \: C_0^{\infty} (\R^{n-1}) \longrightarrow C^{\infty}(\overline{\R}_{+}^n),\: j = 0, 1, 2, 3,$ by
\begin{equation} \label{eq:3.1}
\Pi_j u(x) = (2 \pi h)^{1 -n}\int_{\R^{n-1}}\int_{\gamma} e^{i \langle x , \xi \rangle /h} R_1(x, \xi; z, h) \xi_n^{2j} \hat{u}(\xi')d\xi_n d\xi'/(2 \pi i),\: j = 0, 1,
\end{equation} 
\begin{equation} \label{eq:3.2}
\Pi_j v(x) = (2 \pi h)^{1 -n}\int_{\R^{n-1}}\int_{\gamma} e^{i \langle x , \xi \rangle /h} R_2(x, \xi; z, h) \xi_n^{2j - 4} \hat{v}(\xi')d\xi_n d\xi'/(2 \pi i),\: j = 2, 3,
\end{equation} 
where 
$$\hat{u}(\xi') = \int e^{-i\la x', \xi'\ra/h} u(x') d x'$$ 
is the semi-classical Fourier transform of $u(x').$ Since the poles of the meromorphic functions $\xi_n \to R_j(x, \xi; z, h)$ in the upper half plane are $\lambda_{j, +}, \: \mu_{j, +},\: j = 1,2$, it is easy to see that 
$$\Pi_j u = \frac{h}{i}Op_h(R_1)( u(x') \otimes D_{x_n}^{2j} \delta_{x_n = 0} ) + K_j,\:x_n > 0,\:  j = 0,1$$
with $K_j \equiv 0.$ The equality $M_1 \equiv M_2$ means that the kernel ${\mathcal K}(x, y; z, h)$ of the operator $M_1- M_2$ satisfies
$$|\pa_x^{\alpha} \pa_y^{\beta} {\mathcal K}(x, y; z, h) | \leq C_{\alpha, \beta , N} h^N,\: \forall \alpha,\: \forall \beta, \: \forall N \in \N.$$
We have a similar expression for $\Pi_j u,\: j = 2,3.$ The reader may consult \cite{SjZ} for the properties of these operators.
In particular we have
$$Q_1(x, hD_x) \Pi_j \equiv 0,\: j = 0, 1,\:Q_2(x, hD_x)  \Pi_j \equiv 0, \: j = 2, 3.$$
Here $Q_k(x, hD_x), \: k = 1,2,$ are the operators in the local coordinates $(x', x_n)$ and we use the same notation $q_k(x,\xi)$ for the corresponding principal symbols of $Q_k.$ 

Given $ F =(f_1, f_2 )\in C_0^{\infty} (\overline{\R}_{+}^n: \C^2))$, we define by $\tilde{F} = (\tilde{f_1}, \tilde{f_2})$ the zero extension of $(f_1, f_2)$ to $\R^n$ and set $R(\tilde{F}) =  (R_1(\tilde{f_1}), R_2(\tilde{f_2})).$ \\ 

Introduce the boundary operators
$$B_0 (u, v): = \gamma_0(u - v),$$
$$B_1 (u, v): = \gamma_0h D_{x_n} (u - v),$$
$$B_2 (u, v): = \gamma_0 \Bigl( - (1+m)^2 h^2D_{x_n}^2 u +Y u - h^2D_{x_n}^2 v\Bigr),$$
$$B_3 (u, v): = \gamma_0 \Bigl( -(1+m)^2 h^3D_{x_n}^3 u + Yh D_{x_n} u - h^3D_{x_n}^3 v\Bigr),$$
where the operator $Y(x',  hD_{x'})$ has symbol $Y(x', \xi')$ defined in Section 2.
Let 
$$ B_j(R(\tilde{F})) = W_j, \: j= 0,1,2,3.$$
Setting $W = (W_0, W_1, W_2, W_3),$ we will write the above equalities as ${\mathcal B}(R)(\tilde{F}) = W.$
We will search a parametrix ${\mathcal E}(z)F$ in the form
$${\mathcal E}(z)F = \begin{pmatrix}R_1(\tilde{f_1}) - (\Pi_0 (w_0) + \Pi_1 (w_1)) & & 0\\ 
0 & &R_2(\tilde{f_2}) - (\Pi_2(w_2) + \Pi_3 (w_3)) \end{pmatrix} $$
$$= \begin{pmatrix} R_1(\tilde{f_1}) & & 0 \\ 0 & & R_2(\tilde{f_2}) \end{pmatrix} - \Pi(w)$$
with $w= (w_0, w_1, w_2, w_3)$ and
$$B_j(R(\tilde{F}) - \Pi(w)) \equiv 0,\: j=0,1,2,3.$$
The last relation is equivalent to 
\begin{equation} \label{eq:3.3}
\begin{cases}
B_0\Bigl((\Pi_0(w_0) + \Pi_1(w_1) , (\Pi_2(w_2) + \Pi_3 (w_3)\Bigr) \equiv W_0 ,\\
B_1\Bigl((\Pi_0(w_0)+ \Pi_1(w_1), (\Pi_2(w_2) + \Pi_3 (w_3)\Bigr) \equiv W_1,\\
B_2\Bigl((\Pi_0(w_0) + \Pi_1(w_1), (\Pi_2 (w_2)+ \Pi_3(w_3)\Bigr) \equiv W_2,\\
B_3\Bigl((\Pi_0(w_0) + \Pi_1(w_1), (\Pi_2(w_2) + \Pi_3 (w_3)\Bigr) \equiv W_3.
  \end{cases}
\end{equation}
Here and below for two functions $a$ and $b$ we use the notation $a \equiv b$ if $|a(x')- b(x')| \leq C_N h^N,\: \forall x', \forall N \in N.$ The domain where
this relation holds will be clear from the context.\\

 For simplicity of the notation we write the system (\ref{eq:3.3}) as ${\mathcal M} w \equiv W$ with a matrix-valued h-pseudodifferential operator ${\mathcal M}.$   To construct the parametrix close to the boundary, it is sufficient to prove that the system (\ref{eq:3.3}) is an elliptic one and for this purpose we examine the principal symbol $M$ of ${\mathcal M}.$.  First we treat the case $|\xi'| \geq C_0$, where $C_0 > 0$ is the constant of Proposition 4. Since $|\im z| \geq h^{\ep},$ the roots $\lambda_{j, \pm}, \: \mu_{j, \pm}, \: j = 1,2,$ are simple and we introduce
$$\frac{1}{\pa_{\xi_n} q_1(x', 0, \xi', \lambda_{j, +})} = r_j,\: \frac{1}{\pa_{\xi_n} q_2(x', 0, \xi', \mu_{j, +})} = r_{j+2}, \: j = 1, 2.$$

For simplicity we omit the sign $+$ in the notations below. By applying the theorem of residues, we find 
$$\gamma_0\Bigl((h D_{x_n})^k\Pi_j\Bigr) = Op_h(d_{k,j}), \: j = 0, 1, 2, 3,\: k = 0, 1, 2, 3$$
 with 
$$d_{k, j} = \lambda_1^{2j + k}r_1 + \lambda_2^{2j + k} r_2 + h \tilde{d}_{k, j},\:j = 0, 1, $$
$$ d_{k, j} = \mu_1^{2j - 4 + k}r_3 + \mu_2^{2j - 4 + k} r_4 + h \tilde{d}_{k, j},\: j = 2,3.$$
Here we describe only the principal symbols and $\tilde{d}_{k, j}$ denote lower order symbols. For example the equation 
$$\gamma_0\Bigl[(\Pi_0(w_0) + \Pi_1(w_1)) - (\Pi_2(w_2) + \Pi(w_3))\Bigr] = W_0,$$
modulo lower order terms, has the form
$$\sum_{k=0,1} Op_h(\lambda_1^{2k} r_1 + \lambda_2^{2k} r_2)w_k  - \sum_{k = 2,3}Op_h(\mu_1^{2k-4} r_3 + \mu_2^{2k-4} r_4)w_k = W_0$$ 
and we  find similar expressions for the principal parts of other equations related to the boundary conditions.

Following this argument, the principal symbol $M$ of ${\mathcal M}$ becomes the matrix
$$M(x', \xi', z)  = \begin{pmatrix} r_1 & r_2 & -r_3 & -r_4\\
\lambda_1 r_1 & \lambda_2 r_2 & - \mu_1 r_3 & - \mu_2 r_4\\
(-(1+m)^2 \lambda_1^2 + Y)r_1 & (-(1+m)^2 \lambda_2^2 + Y) r_2 & -\mu_1^2 r_3 & -\mu_2^2r_4\\
(-(1+m)^2 \lambda_1^3 + Y\lambda_1)r_1 & (-(1+m)^2 \lambda_2^3 + Y \lambda_2)r_2 &  -\mu_1^3 r_3 & -\mu_2^3 r_4
\end{pmatrix} \begin{pmatrix} 1 & \lambda_1^2& 0 & 0\\
1 & \lambda_2^2 & 0 & 0\\
0 & 0 & 1 & \mu_1^2\\
0 & 0 & 1 & \mu_2^2 \end{pmatrix}.
$$ 
A simple calculation yields
$$\det M = -(\lambda_1- \lambda_2)^2 (\mu_1-\mu_2)^2 (\lambda_1 +\lambda_2) (\mu_1 +\mu_2) \det \Bigl[\Lambda_1(x', \xi', z) {\rm diag}\: \{r_1, r_2, r_3, r_4\}(x', \xi', z)\Bigr],$$
 where $\Lambda_1(x', \xi', z)$ is the matrix introduced in Section 2. On the other hand,
$$r_1 r_2 = -(\lambda_1 - \lambda_2)^{-2}\prod_{\nu = 1, 2} (\lambda_{\nu, +} - \lambda_{1, -})^{-1} (\lambda_{\nu, +} - \lambda_{2, -})^{-1}$$
$$= -(\lambda_1 - \lambda_2)^{-2} \frac{1}{4 \lambda_1 \lambda_2 (\lambda_1 + \lambda_2)^2},$$
$$r_3 r_4 = - (\mu_1 - \mu_2)^{-2} \prod_{\nu = 1, 2} (\mu_{\nu, +} - \mu_{1, -})^{-1} (\mu_{\nu, +} - \mu_{2, -})^{-1}$$
$$ = -(\mu_1 - \mu_2)^{-2} \frac{1}{4 \mu_1 \mu_2 (\mu_1 + \mu_2)^2}.$$
Thus we obtain
$$\det M = -\frac{\det \Lambda_1(x', \xi', z)}{16\lambda_1 \lambda_2 \mu_1 \mu_2(\lambda_1 + \lambda_2)(\mu_1 + \mu_2)}.$$

According to Proposition 4 for $|\xi'| \geq C_0$, we have $|\det \Lambda_1| \geq C_2$ and since $\lambda_{j} \sim i\sqrt{s},\: \mu_{j,} \sim i\sqrt{s}$, as $|\xi'| \to \infty,$  one obtains
\begin{equation} \label{eq:3.4} 
|\det M (x', \xi', z)| \geq C_4|\xi'|^{-6},\:C_4 > 0, \: |\xi'| \geq C_0.
\end{equation}
Next we denote by ${\mathcal M}$ the h-pseudodifferential  matrix operator with principal symbol $M$ obtained from the equations in (\ref{eq:3.3}). Since the pseudodifferential operator ${\mathcal M}$ is invertible, we can determine  $w_j,\: j = 0,1,2,3,$ from the pseudo-differential system (\ref{eq:3.3}) as $w = {\mathcal M}^{-1} ({\mathcal B}(R(\tilde{F})))$ and this completes the analysis for $|\xi'| \geq C_0$.

For the construction in the region  $z \in Z_0,\:|\im z| \geq h^{\epsilon}, \: |\xi'| \leq C_0$, we need a control of the inverse matrix $M^{-1}(x', \xi', z)$ since we cannot apply Proposition 4. To do this, we follow the argument in \cite{SjZ} with some modifications.\\

Recall that
$$q_1(x, \xi', \xi_n) = (\xi_n^2 + s(x, \xi') - \bar{\omega}_0)(\xi_n^2 + s(x, \xi') - \omega_0),$$
$$ q_2(x, \xi', \xi_n) = \frac{1}{(1 + m(x))^2}(\xi_n^2 + s(x, \xi')
 -(1 + m(x)) \bar{\omega}_0)$$
$$ \times(\xi_n^2 + s(x, \xi') - (1 + m(x))\omega_0).$$
Let $\im z \neq 0$ and let $( q_j(x', 0,\xi', \xi_n) - z)^{-1},\: j = 1, 2,$ be the inverse of the temporal distributions on the whole $\xi_n$ axis. Here we consider $(x', \xi')$ as parameters.
 
 Introduce the distributions 

$$\tilde{u}_j = (q_1(x', 0, \xi', D_{x_n}) - z)^{-1} (D_{x_n}^{2j} \delta), \: j = 0, 1,$$
$$ \tilde{v}_j = (q_2(x', 0, \xi', D_{x_n}) - z)^{-1} (D_{x_n}^{2j - 4} \delta), \: j = 2, 3.$$
Then we have
$$u_j = \tilde{u}_j\vert_{x_n > 0} = \frac{1}{2\pi} \int_{\gamma} e^{i x_n \xi_n} (q_1(x', 0, \xi', \xi_n) - z)^{-1} \xi_n^{2j} d\xi_n , \: j = 0, 1,$$
$$v_j = \tilde{v}_j\vert_{x_n > 0} = \frac{1}{2\pi} \int_{\gamma}  e^{i x_n \xi_n} ( q_2(x', 0, \xi', \xi_n) - z)^{-1} \xi_n^{2j - 4} d\xi_n, \: j = 2,3,$$
where $\gamma = \gamma(x', \xi')$ is the contour used  in the definition of $\Pi_j.$ Set
$$m_{0,0} = u_0(0),\: m_{0, 1} = u_1(0),\: m_{0, 2} = -v_2(0),\: m_{0, 3} = -v_3(0),$$
$$m_{1, 0} = u_1(0),\: m_{1, 1} = (D_{x_n} u_1)(0),\: m_{1,3} = - v_3(0),\: m_{1, 4} = -(D_{x_n}v_3)(0),$$
$$m_{2, k} = -(1+m)^2 (D_{x_n}^2 u_k)(0) + Y u_k(0), k= 0, 1,\: m_{2, j} = -(D_{x_n}^2 v_j)(0),\: j = 2,3,$$
$$m_{3,k} = -(1 + m)^2 (D_{x_n}^3 u_k)(0) + Y (D_{x_n}u_k)(0),\: k = 0, 1, \:m_{3,j} = - (D_{x_n}^3v_j)(0),\: j =2,3.$$
For fixed $(x', \xi')$ we obtain the matrix $M(x', \xi', z) = i\{m_{j,k}(x', \xi', z)\}_{j,k=0,1,2.3}$ introduced above. Given $ a = (a_0, a_1, a_2 , a_3) \in \C^4$, put
$$u = a_0 u_0 + a_1 u_1,\: v = a_2 v_2 + a_3 v_3.$$
Therefore one gets
$$\Bigl(q_1(x', 0, \xi', D_{x_n}) - z\Bigr)u = 0, \Bigl(q_2(x', 0, \xi', D_{x_n}) - z\Bigr) v = 0,\: x_n > 0$$
and
$$u(0) = a_0 m_{0,0} + a_1 m_{0, 1},\: -v(0) = a_2 m_{0, 2} + a_3 m_{0, 3}.$$
$$(D_{x_n}u)(0) = a_0 m_{1, 0}  + a_1 m_{1, 1},\: - (D_{x_n}v)(0) = a_2 m_{1, 2} + a_3 m_{1, 3}.$$
 Similar equalities hold for the third and fourth boundary conditions. Setting $U = (u, v)$, we see that $(u, v)$ satisfies a non-homogeneous boundary problem (\ref{eq:2.17}) with $F_1 = F_2 = 0$ and boundary data $- iMa$ and from (2.18) we deduce the estimate
$$\|U\|_4 \leq \frac{C}{|\im z|} \|M a\|_{\C^4}.$$
On the other hand, 

$$u  = \tilde{u}\vert_{x_n > 0} = ( q_1(x', 0, \xi', D_{x_n}) - z)^{-1} (a_0\delta + a_1 D_{x_n}^2 \delta)\vert_{x_n > 0},$$
$$v =  \tilde{v}\vert_{x_n > 0} = (q_2(x', 0, \xi', D_{x_n}) - z)^{-1} (a_2\delta + a_3 D_{x_n}^2 \delta)\vert_{x_n > 0}$$
with $\tilde{u} = a_0 \tilde{u}_1 + a_1 \tilde{u}_1,\: \tilde{v} = a_2 \tilde{v}_2 + a_3 \tilde{v}_3.$
Since $\tilde{u}, \tilde{v}$ are even, this implies
$$\|a_0 \delta + a_1 D_{x_n}^2 \delta \|_{-4} \leq C_1 \|\tilde{u}\| \leq 2^{1/2} C_1  \|u\| \leq C_2\|u\|_4,$$
$$\|a_2\delta + a_3 D_{x_n}^2 \delta \|_{-4} \leq C_1 \|\tilde{v}\| \leq 2^{1/2} C_1 \|v\|  \leq C_2 \|v\|_4,$$
and we deduce
$$\|U\|_4 \geq C_3\Bigl(\|a_0 \delta + a_1 D_{x_n}^2\delta \|_{-4} + \|a_2 \delta + a_3 D_{x_n}^2 \delta \|_{-4}\Bigr).$$
Since $\|a_0 \delta + a_1 D_{x_n}^2 \delta \|_{-4} + \|a_2 \delta + a_3 D_{x_n}^2 \delta \|_{-4}$ is equivalent to $C_4 \|a\|_{\C^4},$ we conclude that with a constant $C_5 > 0$ which is uniform for $|x'| \leq \rho, |\xi'| \leq C_0$ we have
\begin{equation} \label{eq:3.5}
\|M^{-1}(x', \xi', z) \| \leq C_5/ |\im z|.
\end{equation}
It is clear also that we have
\begin{equation} \label{eq:3.6}
\|M(x', \xi', z)\| \leq C_5 / |\im z|
\end{equation}
uniformly with respect to $|x'| \leq \rho,\: |\xi'| \leq C_0.$

Now write $M(x', \xi', z) = (M_{i, j})_{i, j = 1,2, 3 , 4}$ for the principal symbol of ${\mathcal M}$. To study the symbols $(M_{i, j})$,  note that
$$r_{1}\lambda_{1, +}^k = \lambda_{1, +}^k\Bigl((\lambda_{1, +} - \lambda_{1, -})(\lambda_{1, +} -\lambda_{2, +})(\lambda_{1, +} - \lambda_{2, -})\Bigr)^{-1}$$   
$$= \lambda_{1, +}^{k-1}\frac{(\lambda_{1, +} + \lambda_{2, +}) (\lambda_{1, +} + \lambda_{2, -})}{8 (r + i \im z)},$$
with $r = \re z - (\im \om)^2$ and we have  similar expressions for $r_2\lambda_{2, +}^k, r_3 \mu_{1, +}^k, r_4 \mu_{2, +}^k.$
On the other hand, it is easy to see that
$$\lambda_{j, \pm} \in S_{\epsilon}^{0, 1}, \mu_{j, \pm} \in S_{\epsilon}^{0, 1}, j = 1,2.$$
In fact, for $|\xi'| \gg 1$ we have
$$|\partial_{x'}^{\alpha} \partial_{\xi'}^{\beta}  \lambda_{j, \pm}| \leq C_{\alpha, \beta} ( 1+ |\xi'|)^{1 - |\beta|},\: \forall \alpha, \forall \beta,$$
while for bounded $|\xi'|$ we have
$$|\partial_{x'}^{\alpha} \partial_{\xi'}^{\beta}  \lambda_{j, \pm}| \leq C_{\alpha, \beta} |\im z|^{-|\alpha|- |\beta|} \leq C_{\alpha, \beta}'h^{-\epsilon(|\alpha| + |\beta|)},\: |\alpha|+  |\beta| \geq 1.$$
We have similar estimates for $\mu_{j, \pm}.$
The reader may consult also Section 2 in \cite{V} for the estimates of the symbols $\lambda_{j, \pm}, \mu_{j, \pm}.$ Concerning the calculus of $h$-pseudo-differential operators in the class $S_{\epsilon}^{k, m}$, recall that if $a \in S_{\epsilon}^{k_1, m_1},\: b \in S_{\epsilon}^{k_2, m_2}$ we have $a b \in S_{\epsilon}^{k_1 + k_2, m_1 + m_2}$ and the principal symbol of $a \#_h b$ is $ab$
(see for more details Appendix A in \cite{HK3}). Thus  we conclude that $r_j \lambda_{j, +}^k \in S_{\epsilon}^{\epsilon, k + 1  }, r_{j + 2} \mu_{j, +}^k \in S_{\epsilon}^{\epsilon, k + 1},\: j = 1, 2.$ Now it is easy to see that
$$M_{1, j} \in S_{\epsilon}^{\epsilon, 1}, j = 1, 3, M_{1, j} \in S_{\epsilon}^{\epsilon, 3}, j = 2, 4$$
and similarly we get $M_{i,j} \in S_{\epsilon}^{\epsilon, m(i, j)}$. The  precise orders $m(i, j)$ are not important when we consider these symbols for bounded $|\xi'|$. Moreover, the matrix $|\Im z| M(x', \xi', z)$ is bounded and this agrees with the estimate (\ref{eq:3.6}). Next for the inverse matrix $M^{-1} = (s_{i, j})_{i, j = 1,2,3,4}$ we can obtain a similar description of $s_{i, j}$ and we omit the details since this is not necessary for our argument. Indeed, for the calculation of the trace in Section 4 (see Lemma 6) we need to know the behavior of the symbols for bounded $|\xi'|.$ Next we are interested on the trace of the kernels and not on their derivatives. \\

 Set ${\mathcal B}_k=B_k(R),\: k = 0,1,2,3$ and as above let ${\mathcal M}$ be the matrix h-pseudo-differential operator with principal symbol $M$. Then one determines
 $$w = (w_1,...,w_4) = {\mathcal M}^{-1} \Bigl({\mathcal B}_0(\tilde{F}), {\mathcal B}_1(\tilde{F}), {\mathcal B}_2(\tilde{F}), {\mathcal B}_3(\tilde{F})\Bigr) = {\mathcal M}^{-1} {\mathcal B}(\tilde{F}).$$
 This completes the construction of the parametrix ${\mathcal E}(z)$ for $x$ close to the boundary in the case $|\xi'| \leq C_0,\: |\im z| \geq h^{\epsilon}.$\\

To construct a global parametrix, choose a finite number of points $x_k\in \bar {\Omega},\: k =1,...,L,$ so that $x_k$ with  $k = 1,...,L',$ lie in the interior of $\Omega$, while $x_k$ with $k = L'+ 1,...,L,$ belong to $\pa\Omega.$ Let $U_k$ be open neighborhoods of $x_k.$ Assume that we have a covering $\cup_{k =1}^L U_k$ of $\bar\Omega$ and suppose that in every $U_k,\: k= L'+1,...,L,$ we can introduce normal coordinates. Next, let
$\Phi_k \in C_0^{\infty}(U_k),\: k = 1,..., L,$ form a partition of unity on $\bar{\Omega}.$ Introduce $\Psi_k \in C_0^{\infty} (U_k),\: k = 1,...L'$ with $\Psi_k = 1$ near supp $\Phi_k.$ We may choose $\Psi_k(x),\: k = 1,...,L',$ so that supp $\Psi_k \cap \pa \Omega = \emptyset.$ Next let $\Psi \in C_0^{\infty}(\R^n)$ be a function which is equal to 1 on $\cup_{j= 1}^L {\rm supp}\: \Phi_k$ and
let $ {\mathcal M}^{-1} {\mathcal B}(\Phi_k \tilde{F})),\: k = L'+ 1,...,L $ be determined as above (Here ${\mathcal M}^{-1}, {\mathcal B}$ depend of $k$ but we omit this in the notation).

We define a global parametrix by
$${\mathcal E}(z) \tilde{F} = \sum_{k=1}^{L'}\Psi_k R (\Phi_k \tilde{F}) +  \sum_{k= L'+ 1}^L\Psi \Bigl( R(\Phi_k \tilde{F}) - \Pi({\mathcal M}^{-1} ({\mathcal B}(\Phi_k \tilde{F}))\Bigl).$$
We get
$$(Q - z){\mathcal E}(z) \equiv I,$$
$$B_j ({\mathcal E}(z)) \equiv 0,\: j = 0, 1, 2, 3.$$
By a modification of ${\mathcal E}(z)$ by an operator $\equiv 0$, we can arrange the boundary conditions
$$B_j({\mathcal E}(z)) = 0, \: j= 0,1,2,3,$$
preserving the first relation which means that $(Q - z) {\mathcal E}(z) - I$ has a kernel $K(x, y; z, h)$ such that
$$|\pa_x^{\alpha} \pa_y^{\beta} K(x, y; z, h)| \leq C_{\alpha, \beta, N} h^N,\:  (x, y) \in (\Omega \times \Omega), \: \forall \alpha, \forall\beta, \forall N.$$
To do this we use the analysis of the non-homogeneous problem (2.17) and the functions $H_j(x_n), j = 0, 1, 2, 3,$ introduced in Section 2. If we have a solution of (2.17) with operators $f_j(x') \equiv 0,\:j = 0,1,2,3,$ we can reduce the problem to one with  homogeneous boundary conditions changing the right hand part $(F_1, F_2)$ by a negligible term $G \equiv 0.$

Now we pass to the construction of a microlocal parametrix in the region (A). Here we must work in a region including an interval lying on the real axis, so the estimates with ${\mathcal O}(|\im z|^{-1})$ are not useful. The advantage is that for $|\xi'| \geq C_0$, where $C_0 > 0$ is the constant of Propositions 4 and 5,  we have uniform lower bounds for $|\det\Lambda_1|,\:|\det\Lambda_2|$. More precisely, for $k = 1,...,L',$ we will take $|\xi| \geq C_1 > 0$, while for $k = L'+ 1,...,L$ we choose  $|\xi'| \geq C_0$.  Moreover, we will show that this microlocal parametrix depends {\it holomorphically} of $z \in Z_0.$ Assume that $z \in Z_0,\: \re z \leq (\im \omega_0 + \delta)^2,$ where $\delta$ and $\omega_0$ satisfy the conditions mentioned in the beginning of the section.  Choosing a  big constant $C_1 > 0$ (depending on $Z_0$), let $\psi(\xi) \in C_0^{\infty}(\R^n)$ be such that $\psi(\xi) = 1$ for $|\xi| \leq C_1.$  
The construction of $R( (1- \psi(hD_x))\Phi_k(x)\tilde{F}),\: k = 1,...,L',$ so that
\begin{equation} \label{eq:3.10}
(Q- z) \Psi_kR( (1 - \psi(h D_x)) \Phi_k) \equiv (1 - \psi(h D_x)) \Phi_k ,\: k = 1,...,L'
\end{equation}
is standard and it is clear that $R((1 - \psi(hD_x))\Phi_k(x))$ depends holomorphically of $z \in Z_0$. Next,  for $k = L'+ 1,...,L,$ we 
choose a fixed function $\chi(\xi') \in C_0^{\infty} (\R^{n-1})$ such that $\chi(\xi') = 1$ for $|\xi'| \leq C_0$ and we construct 
$R$ so that 

$$(Q- z) \Bigl[\Psi R( (1 - \chi(h D_{x'})) \Phi_k )\Bigr] \equiv (1 - \psi(h D_{x'})) \Phi_k,\: k = L'+ 1,...,L.$$

Obviously, the terms 
$$B_j\Bigl( R ((1 - \chi(h D_{x'}) \Phi_k)\Bigl),\: j = 0,1,2,3$$
are holomorphic with respect to $z$. To treat the boundary terms
 $\sum_{j = 0, 1}\Pi_j (w_j)$ and $\sum_{j= 2, 3} \Pi_j (w_j)$, we must know that the inverse matrix $M^{-1}(x', \xi', z)$ is holomorphic in $z$ for $|\xi'| \geq C_0$.  As in \cite{SjZ}, for $|\xi'| \geq C_0$ we choose $\gamma(x', \xi') \subset \Omega_1(x', \xi', z) \cap \{|\xi'| \geq C_0,\:\im \xi \geq \eta \la \xi'\ra \}$ with $\eta > 0$ small enough. Then for $|\xi'| \geq C_0$ the components $M_{j, k}(x', \xi', z)$ of the matrix $M(x', \xi', z)$ are holomoprphic with respect to $z \in Z_0$ and so is $\det M(x', \xi',z).$  In the case $z \neq (\im \om)^2$, applying Proposition 4 for $\det \Lambda_1(x', \xi', z)$, we obtain (\ref{eq:3.4}). Since $\det M(x', \xi', z)$ is continuous with respect to $z$, the estimate (\ref{eq:3.4}) remains true for $z = (\im \om)^2$ and $|\xi'| \geq C_0.$ Another proof of this statement can be obtained by applying Proposition 5 combined with an expression of  $M(x', \xi', z)$ involving the double roots $\lambda_{+},\: \mu_{+}.$

The matrix $M(x', \xi', z)$ is invertible and holomorphic in $z$ for $z \in Z_0$ and $|\xi'| \geq C_0$, hence we obtain  immediately that $M^{-1}(x', \xi', z)$ is holomorphic with respect to $z \in Z_0$ for  $|\xi'| \geq C_0$. Following this way we deduce that the matrix-valued operator ${\mathcal M(z)}^{-1}$ with principal symbol $M^{-1}(x', \xi', z)$ defined above composed with  ${\mathcal B}((1- \chi(h D_{x'}))\Phi_k(x))= \{B_j(R(1 - \chi(h D_{x'})\Phi_k)\}_{j=0,1,2,3}$ yields an operator holomorphic with respect to $z.$ Finally, we obtain a microlocal parametrix
\begin{equation} \label{eq:3.11}  
(Q- z)E_k(z) \equiv (1 - \chi(h D_{x'})) \Phi_k(x),\:  B_j(E_k(z)) = 0,\:k = L'+ 1,...,L,\: j = 0,1,2,3
\end{equation}
depending holomorphiically of $z \in Z_0.$

\section{The trace of $f(Q)$}

In this section we use the notations of the previous one and we assume (\ref{eq:cond}) fulfilled. Let $f \in C_0^{\infty}( ]-\infty, (\im \omega_0 + \delta)^2[)$ be a cut-off function with sufficiently small $\delta > 0.$ Let $z$ belong to a fixed compact set  $Z_0$ and let $\re z \leq (\im \omega_0 + \delta)^2.$ We choose an almost analytic extension $\tilde{f}(z) \in C_0^{\infty} (\{ z \in Z_0:\: \re z < (\im \omega_0 + \delta)^2\})$ of $f$ so that $\bar{\pa}\tilde{f}(z) = \Oc(|\im z|^N),\: \forall N \in \N.$  Next, as in Section 4 of \cite{SjZ}, we use the formula
\begin{equation} \label {eq:4.1}
f(Q) = - \frac{1}{\pi} \int \bar{\pa} \tilde{f}(z) (z - Q)^{-1} L(dz),
\end{equation}
where $L(dz)$ is the Lebesgue measure in $\C.$ Notice that $(z - Q)^{-1}$ is a matrix-valued operator and so is the operator $f(Q)$.\\

Following our construction in Section 3, the global parametrix is a sum of $L$ terms.  For $k = 1,..., L',$ from (\ref{eq:3.10}) we get
$$-\Psi_k R((1 - \psi(hD_x)) \Phi_k) = (z- Q)^{-1} (1 - \psi(hD_x))\Phi_k + K_k(z)$$
with an operator $K_k(z)$ having trace class norm $[K_k(z)] = {\mathcal O}(|\im z|^{-1}h^{\infty}).$ The operator on the left hand side is 
holomorphic with respect to $z$ so the contribution of $(z- Q)^{-1}(1 - \psi(hD_x))\Phi_k$ in (\ref{eq:4.1}) is negligible. For $k = L'+ 1,...,L,$ we apply the same argument exploiting now (\ref{eq:3.11}) and the fact that $E_k(z)$ is holomorphic with respect to $z$. Thus we reduce the analysis to the examination of the sum
$$f(Q) = -\frac{1}{\pi} \sum_{k= 1}^{L'} \int\bar{\pa}\tilde{f}(z) (z - Q)^{-1}\psi(h D_{x})\Phi_k(x) L(dz)$$
$$ -\frac{1}{\pi} \sum_{k= L'+1}^{L} \int\bar{\pa}\tilde{f}(z) (z - Q)^{-1}\chi(h D_{x'})\Phi_k(x) L(dz) + K_L$$
with $[K_L] = {\mathcal O}(h^{\infty})$ and $\psi(\xi),\: \chi(\xi')$ introduced in Section 3.\\

Before going to a second reduction, it is necessary to study the traces of some operators. We start with a description of the distribution kernel $\Pi_j(x, y'; z, h)$ of the operator $\Pi_j$ given by the following
\begin{lem} Fix $0<\epsilon<1/2$ small enough. Let $\vert\Im z\vert\geq h^\epsilon$ and let $z$ belong to a compact set $Z_0$. Modulo a term having order ${\mathcal O}(h^\infty)$ in the trace class norm, we have the following representation  
$$\Pi_j(x, y' ;z, h) = (2\pi h)^{-(n-1)} \int e^{i (x'- y')\xi'/h} z_j(x, \xi'; z, h) d\xi'.$$
The symbol $z_{j}$ satisfies the estimates
\begin{eqnarray} 
|\im z| |\pa_{x'}^{\alpha'}  \pa_{\xi'}^{\gamma'} (h D_{x_n})^k z_{j}(x, \xi'; z, h)| \nonumber \\
\leq C_{\alpha', \gamma', k} \la \xi'\ra ^{ 2j + 1+ k - |\gamma'|} |\im z|^{- |\alpha'| - |\gamma'|}\Bigl(e^{- x_n \Im  \lambda_1 /h}+e^{- x_n \Im \lambda_2 /h}\Bigr),\: j = 0, 1 \label{eq:4.2} 
\end{eqnarray}
and
\begin{eqnarray} 
|\im z| |\pa_{x'}^{\alpha'} \pa_{\xi'}^{\gamma'} (h D_{x_n})^k z_{j}(x, \xi'; z, h)| \nonumber \\
\leq C_{\alpha', \gamma', k} \la \xi'\ra ^{ 2j -3 + k - |\gamma'|} |\im z|^{-|\alpha'| - |\gamma'|}\Bigl(e^{- x_n \Im  \mu_1 /h}+e^{- x_n \Im \mu_2 /h}\Bigr), j = 2, 3. \label{eq:4.3} 
\end{eqnarray}
Moreover, for $|\xi'| \geq C_1$ (depending on the compact set $Z_0$) the above estimates  hold without $|\im z|^{-|\alpha'| - |\gamma'|}$ in the right hand part. 
\end{lem}

{\bf Proof.} According to (\ref{eq:3.1}) and (\ref{eq:3.2}), for the principal symbols $z_{p, j}$ we have the representation
$$z_{p,j}(x,\xi',z, h)=\int_{\gamma} e^{i \langle x_n , \xi_n \rangle /h} q_1(x, \xi, z) \xi_n^{2j} d\xi_n /(2 \pi i),\: j =0,1,$$
$$z_{p,j}(x,\xi',z, h)=\int_{\gamma} e^{i \langle x_n , \xi_n \rangle /h} q_2(x, \xi, z) \xi_n^{2j-4} d\xi_n /(2 \pi i),\: j =2,3.$$
We treat only the case $j = 0, 1$. Applying the theorem of the residues, we get
\begin{equation} \label{eq:4.4}
z_{p,j}(x,\xi'; z, h)=\sum_{\nu = 1, 2} \lambda_{\nu, +}^{2j}r_{\nu} e^{ix_n\lambda_{\nu, +}/h},\: j = 0, 1,
\end{equation}
where $r_{\nu} = \Bigl(\pa_{\xi_n} q_1 (x', x_n, \xi', \lambda_{\nu, +})\Bigr)^{-1}.$
We established in Section 3 the estimates for $\partial_{x'}^{\alpha'} \partial_{\xi'}^{\gamma'} (r_{\nu}\lambda_{\nu, +}^k)$ which imply the estimates for the derivatives $\partial_{x'}^{\alpha'} \partial_{\xi'}^{\gamma'}z_{p, j}$.  On the other hand,  when we take the derivatives $(h D_{x_n})^k \lambda_{\nu, +}$, we have the factor $\Bigl(\frac{h}{|\im z|}\Bigr)^k$ which is bounded for $|\im z| \geq h^{\epsilon}.$

We remark that a similar estimates are true for the full symbol of the kernels. In fact,  modulo a term having order ${\mathcal O}(h^\infty)$, $R_1(x,\xi; z,h)$ is a sum of terms of the form $h^N b_N(x,\xi ;z) (q_1-z)^{-2N-1}$, where $b_N(x,\xi; z)$ a polynomial in $z$ with smooth coefficients (see Proposition 8.6 and Theorem 8.7 in  \cite{DS}). Therefore an application of the theorem of residues leads to a sum of terms involving derivatives of order $2N$. We treat below the behavior for bounded $|\xi'|$ when the factor $|\Im z|^{-1}$ is important. Consider a typical term
$$h^N \frac{\partial^{2N}}{\partial\xi_n^{2N}}\Bigl[ \frac{\xi_n^{2j} b_N(x, \xi; z)}{\Bigl(\xi_n - \lambda_{1, -})(\xi_n- \lambda_{2, +})(\xi_n - \lambda_{2, -})\Bigr)^{2N + 1}}\Bigr]\Bigl\vert_{\xi_n = \lambda_{1, +}}.$$
The product in the dominator yields an expression having the form
$$\prod_{\alpha + \beta_1 + \beta_2 \leq 4N + 1} B_{\alpha, \beta_1, \beta_2} (\lambda_{1, +} - \lambda_{1, -})^{-\alpha} (\lambda_{1, +} - \lambda_{2, +})^{-\beta_1} (\lambda_{1, +} - \lambda_{2, -})^{-\beta_2}$$
$$= \prod_{\alpha + \beta_1 + \beta_2 \leq 4N + 1} B_{\alpha, \beta_1, \beta_2} (2\lambda_{1, +})^{-\alpha} (\lambda_{1, +}+\lambda_{2, +})^{\beta_1} (\lambda_{1, +} + \lambda_{2, -})^{\beta_2} \Bigl(r + i \Im z\Bigr)^{-\frac{1}{2} (\beta_1+\beta_2)}. $$
Next, $\frac{1}{2}(\alpha + \beta_1 + \beta_2) \leq 2N + 1$ and
$$h^N |\Im z|^{-\frac{\alpha}{2}}|(r + i \Im z)^{-\frac{1}{2}(\beta_1 + \beta_2)}| \leq h^{N}|\Im z|^{-2N - 1}\leq C h^{N(1 - 2 \epsilon)}|\Im z|^{-1}$$ 
since $\vert\Im z\vert^{-2} \leq h^{- 2\epsilon}$. Thus for $\alpha' = \gamma'= 0$ and bounded $|\xi'|$ we obtain  (\ref{eq:4.2}) for the full symbol of the kernel. The analysis of the derivatives $\partial_{x'}^{\alpha'}\partial_{\xi'}^{\gamma'}$ follows the same argument.  The proof of (\ref{eq:4.3}) is similar. \hfill\qed\\

Now we pass to the analysis of the kernels of the operators

$$Y_{k, j} = \Bigl((hD_{x_n})^ k R_j(x,hD_x; z,h)\Bigr)\bigg\vert_{x_n=0},\: k= 0, 1, 2, 3,\: j = 1, 2.$$ 

The kernels of these operators have the form
$$(2\pi h)^{-n}\int_{{\R}^{n-1}_{\xi'}}e^{i(x'-y')\xi'/h}\int_{\mathbb R} e^{-iy_n\xi_n/h}\xi_n^k R_j(x',0,\xi',\xi_n ;z,h)  d\xi_n d \xi'.
$$
For simplicity we treat below the case $j = 1$. In Section 3 it was shown that
 $$R_1(x',0,\xi',\xi_n; z, h)= \sum_{\nu = 0}^N h^{\nu} c_{\nu}(x',0, \xi',\xi_n; z)+h^{N + 1}c_{N + 1}(x',0, \xi',\xi_n;z),$$ 
with $c_{N + 1}\in S_{\epsilon}^{\epsilon(N + 2),-4-N - 1} $ and $c_{\nu}$ are holomorphic in the domain $\Omega(x, \xi', z)$ defined in Section 3. Clearly,
$$(2\pi h)^{-n}\int_{ {\R}^{n-1}_{\xi'}}e^{i(x'-y')\xi'/h}\int_{\mathbb R} e^{-iy_n\xi_n/h}\xi_n^k  h^{N+1}c_{N +1}(x',0, \xi',\xi_n; z)  d\xi_n d\xi'={\mathcal O}(h^{M(N)}),
$$
uniformly for $\vert \im z\vert >h^\epsilon$  with $M(N) \to \infty$ as $N \to \infty.$
Let $\Gamma_\rho(x',\xi')=[-\rho, \rho]\cup \gamma_\rho(x',\xi')$ be a closed, simple loop, which encircles $\lambda_{j,-}, \mu_{j,-}$ in the positive sense.
Here $ \gamma_\rho(x',\xi')$ is a curve included in $\{\xi_n\in{\C};\: \im \xi_n\leq 0, \vert \xi_n-\lambda_{j,-}\vert \geq \eta  \la \xi'\ra, \:\vert \xi_n-\mu_{j,-}\vert \geq \eta  \la \xi'\ra,\: j = 1,2\},\: \eta > 0 $, and $\rho > R_0(x',\xi')$ is a large positive constant.
We fix $(x', \xi')$ and recall that $y_n > 0.$ We apply the theorem of the residues to the integral
$$G_{\nu}(x', y_n, \xi';z, h) = \frac{1}{2 \pi i}\int_{\Gamma_{\rho}} e^{-i y_n \xi_n/h} \xi_n^k c_{\nu}(x', 0, \xi', \xi_n; z) d \xi_n.$$
For $\nu = 0$ we have $c_0 = (q_1 - z)^{-1}$ and we obtain 
$$\sum_{j = 1,2} \lambda_{j, -}^k r_{j, -} e^{-i y_n \lambda_{j, -}/h}$$
with $r_{j, -} = \Bigl(\partial_{\xi_n} q_1(x', 0, \xi', \lambda_{j, -})\Bigr)^{-1}.$ The analysis  of the derivatives of this symbol is
completely similar to that of (\ref{eq:4.4}) and we deduce
$$|\im z| |\partial_{x'}^{\alpha'} \partial_{\xi'}^{\gamma'}(hD_{y_n})^m (G_{0})|\leq C_{\alpha', \gamma', m, j} \langle \xi'\rangle^{1 + m + k- |\gamma'|} |\im z|^{-|\alpha'|- |\gamma'|} e^{-y_n c |\Im z| /h}.$$
We can treat the lower order symbols $c_{\nu}$ by a similar argument. Next the integral over $\gamma_{\rho}$ goes to 0 as $\rho \to +\infty$ and taking the limit $\rho \to +\infty$, we obtain an estimate for the kernels of $Y_{k, j}$. Consequently, for the  kernels of ${\mathcal B}_k$ we have the following
\begin{lem} Under the assumptions of Lemma $4$ the kernel of the operator ${\mathcal B}_k = B_k(R),\: k = 0, 1, 2, 3$ have the form
$${\mathcal B}_k(x', y; z, h) = (2\pi h)^{-n} \int e^{i \langle x'- y', \xi'\rangle/h} b_k(x', y, \xi'; z) d\xi',$$
with
\begin{equation} \label{eq:4.5}
|\im z| |\pa_{x'}^{\alpha'} \pa_{y}^{\beta} \pa_{\xi'}^{\gamma'} (h D_{y_n})^m b_k(x', y, \xi';z)| 
\leq C_{\alpha', \beta, \gamma', j} \la \xi'\ra ^{1 + m + k - |\gamma'|} |\im z|^{-|\alpha'| -|\gamma'|} e^{- y_n c|\im z| /h}.
\end{equation}
Moreover, for $|\xi'| \geq C_1$ the estimates hold without the factor $|\im z|^{-|\alpha'| - |\gamma'|}$.
\end{lem}
Next we need to examine the trace of the boundary terms $\Pi_k(w_k), k = 0,1,2,3$ involved in the parametrix ${\mathcal E}(z).$ Using the notation of the previous section, recall that
$$w = {\mathcal M}^{-1} \Bigl(\{{\mathcal B}_k(R)\}_{k = 0,1,2,3}\Bigr)= {\mathcal M}^{-1} {\mathcal B}(\tilde{F}),$$
so to obtain the kernel of the vector-valued operator ${\mathcal M}^{-1} {\mathcal B} := \{A_0, A_1, A_2, A_3\}$ we must compose the kernels of the matrix-valued h-pseudodifferential operator ${\mathcal M}^{-1}$ and the vector-valued operator ${\mathcal B}.$ Let $\chi_1(x,\xi') \in C^\infty_0({\mathbb R}^n\times {\mathbb R}^{n-1})$. For our argument we need to know the traces of the operators $(\Pi_j A_j + \Pi_{j + 1} A_{j+1})\chi_1(x, h D_{x'}),\: j = 0, 2$ for $z \in Z_0$ and $|\im z | \geq h^{\epsilon}.$ To do this, we must compose the kernels of the operators $\Pi_j, {\mathcal M}^{-1}$ and ${\mathcal B}$. For bounded $|\xi'|$ taking into account Lemma 4 and Lemma 5, we conclude  that in the derivatives of the 
kernels the terms $|\im z|^{-k}$ can be estimated by $h^{-k \epsilon}$, while the terms $\langle \xi'\rangle^{-m}$ will be bounded. More precisely we have the following

\begin{lem}
Let $\chi_1(x,\xi') \in C^\infty_0({\mathbb R}^n\times {\mathbb R}^{n-1})$. Under the assumptions of Lemma $4$, the operators 
$$(\Pi_j A_j + \Pi_{j + 1} A_{j+1})\chi_1(x, h D_{x'}),\: j = 0, 2,$$
  are trace class and
\begin{equation} \label{eq:4.6}
{\rm tr}\Bigl((\Pi_j A_j + \Pi_{j+1} A_{j+1})\chi_1(x, h D_{x'})\Bigr)={\mathcal O}(h^{1-n-4\epsilon}),\,\,j=0,2.
\end{equation}
\end{lem}

{\bf Proof.} Let $\ell_j(x, y; z, h)$ be the kernel of   $(\Pi_j A_j + \Pi_{j+1} A_{j+1})\chi_1(x, h D_{x'}),\: j = 0, 2$. Since all terms of  the operator
 $(\Pi_j A_j + \Pi_{j+1} A_{j+1})\chi_1(x, h D_{x'})$ are $h$-pseudodifferential operators in $x'$, we deduce from (\ref{eq:4.2}),
(\ref{eq:4.3}) and (\ref{eq:4.5}) that modulo a term of order ${\mathcal O}(h^\infty)$ in trace class norm we have
$$\ell_j(x,y;z, h) = (2\pi h)^{-n} \int e^{i\la x'- y', \xi'\ra/h} c_j(x, y, \xi'; z)\chi_1(y, \xi') d\xi'$$
with
$$|\pa_{x'}^{\alpha'} \pa_{y'}^{\beta'} \pa_{\xi'}^{\gamma'} (h D_{x_n})^k (hD_{y_n})^m c_j(x, y, \xi'; z)|\leq C_{\alpha', \beta', \gamma', k, m}  h^{-3\epsilon- |\alpha'| - |\gamma'|} \exp(-(x_n + y_n)h^\epsilon/Ch).$$ 
Here for $\alpha'= \gamma' = 0$ we have twice the factor $|\im z|^{-1}$ from (\ref{eq:4.2}) and (\ref{eq:4.5}) and once more the same factor  from the estimates of ${\mathcal M}^{-1}.$
This explains the factor $h^{-3 \epsilon}.$ To calculate the trace we must take the trace $\ell(x, x; z, h)$ and integrate. 
This yields  (\ref{eq:4.6}). \hfill\qed\\

In the case when $z = z_0$ is fixed with $\im z_0 > 0$ the estimates in Lemma 4 and Lemma 5 imply an estimate for the kernel

 $\ell_j(x, y; z, h)$ having the form
$$
|\pa_{x'}^{\alpha'} \pa_{y'}^{\beta'} \pa_{\xi'}^{\gamma'} (h D_{x_n})^k )(hD_{y_n})^m \ell_j(x, y; z, h))| \leq C_{\alpha', \beta', \gamma', k, m}h^{-n}\exp(-(x_n + y_n)/Ch) + {\mathcal O}(h^{\infty}).
$$
with constants depending of $\im z_0.$
Consequently, in this case we have

\begin{equation} \label{eq:4.7}
\Bigl\|(\Pi_j A_j + \Pi_{j+1} A_{j+1})\chi_1(x, h D_{x'})\Bigr\|_{tr}={\mathcal O}(h^{1-n}),\,\,j=0,2.
\end{equation}

To do a second reduction, consider a function $\chi_k(x, \xi')$ with support in $\{(x, \xi'):\: x \in U_k, \: |\xi'| \leq C_0\}$, $U_k$ being a small neighborhood of $x_k \in \pa \Omega.$  We fix $z_0$ with $ \im z_0 > 0, \: \re z_0 < (\im \omega_0)^2$ and we construct a parametrix $E(z_0)$ for the operator $(Q- z_0)$ as in Section 3. Thus we obtain
\begin{equation} \label{eq:4.8}
(z_0 - Q)^{-1} \chi_k(x, h D_{x'}) =  -E(z_0) \chi_{k}(x, h D_{x'}) + K_k
\end{equation}
with trace class norm $[K_k] = {\mathcal O}(h^{\infty}).$\\ 

Next, we are going to repeat the arguments in \cite{SjZ} with some minor modifications. The analysis of $R(x, hD_x; z_0, h) \chi_k(x, hD_{x'})$ is easy and we obtain a trace class operator with trace class norm ${\mathcal O}(h^{-n}).$  On the other hand, for fixed $z_0$ with $\im z_0 > 0$ we apply  (\ref{eq:4.7}). 
Thus for $\im z \neq 0$ and $z \in Z_0$  we have
$$(z - Q)^{-1} \chi_k(x, hD_{x'}) = (I + (z_0 - z) (z - Q)^{-1}) (z_0 - Q)^{-1} \chi_k(x, hD_{x'})$$
 and taking into account the estimate $\|(Q - z)^{-1}\| \leq \frac{1}{|\im z|}$, we get
\begin{equation} \label{eq:4.9}
\Vert (Q- z)^{-1}\chi_k(x, h D_{x'})\Vert_{\rm tr} \leq C | \im z|^{-1} h^{-n }.
\end{equation}
On the other hand, for $\tilde{f}(z)$ we have the estimate $\bar{\pa}\tilde{f}(z) = {\mathcal O}(|\im z|^N),\: \forall N.$ Combining this with (\ref{eq:4.9}),  leads to
\begin{equation} \label{eq:4.10}
 \int  \bar{\pa} \tilde{f}(z) (z - Q)^{-1}\chi_k(x, hD_{x'}) L(dz)  =  \int_{|\im z| \geq h^{\epsilon}}  \bar{\pa} \tilde{f}(z) (z - Q)^{-1}\chi_k(x, hD_{x'}) L(dz) + K_{\epsilon},
\end{equation}
with $[K_{\epsilon}] = {\mathcal O} (h^{\infty})$ since the trace of the operator corresponding to the integral over $|\im z| \leq h^{\epsilon}$ is negligible.\\

Now we apply the partition of unity $\{ \Phi_k(x)\}_{k=1,...,L}$ introduced in Section 3. Set
$$\chi_k(x, \xi) = \psi(\xi) \Phi_k(x), \: k = 1,..., L', \: \chi_k(x, \xi') = \chi(\xi') \Phi_k(x),\: k = L'+ 1,...,L$$
with the functions $\psi(\xi),\: \chi(\xi')$ introduced in Section 3. 
Let $E_k(z)$ be the parametrix of $(Q- z)$ constructed in Section 3 for $z$ in a compact set $Z_0$ and $|\im z| \geq h^{\epsilon}$ with $\epsilon > 0$ small enough for which we have 
\begin{equation} \label{eq:4.11}
(z-Q)^{-1}\chi_k(x, h D_{x'}) = -E_k(z)\chi_k(x, hD_{x'}) +K_k(z),\: k = L'+ 1,...,L
\end{equation}
 with $[K_k(z)]={\mathcal O}(h^\infty)$, uniformly on $z$. For $k = 1,...,L',$ we have a similar expression and the paramterix $E_k(z)$ has no boundary terms.
 
According to (\ref{eq:4.10}),  modulo terms with trace class norm ${\mathcal O}(h^{\infty})$,  we get 
\begin{eqnarray}
{\rm tr}(f(Q))= {\rm tr}\: \Bigl[ - \frac{1}{\pi}\sum_{\nu=1}^2\sum_{k = 1}^{L'} \int_{|\im z| \geq h^{\epsilon}}  \bar{\pa} \tilde{f}(z) R_{\nu, k}(x, hD_x,; z, h)\chi_k(x, hD_{x}) L(dz)\Bigr] \nonumber\\
+{\rm tr}\: \Bigl[ - \frac{1}{\pi}\sum_{\nu=1}^2\sum_{k = L'+ 1}^{L} \int_{|\im z| \geq h^{\epsilon}}  \bar{\pa} \tilde{f}(z) R_{\nu, k}(x, hD_x,; z, h)\chi_k(x, hD_{x'}) L(dz)\Bigr]\nonumber\\
-{\rm tr}\: \Bigl[  \frac{1}{\pi}\sum_{\nu=0}^2\sum_{k = L'+ 1}^{L} \int_{|\im z| \geq h^{\epsilon}}  \bar{\pa} \tilde{f}(z) \Bigl(\Pi_{\nu} A_{\nu, k} + \Pi_{\nu + 1} A_{\nu + 1, k}\Bigr)\chi_k(x, h D_{x'}) L(dz)\Bigr] : =(I)+(II)+ (III).\nonumber
\end{eqnarray}
Here $R_{\nu, k} = \Psi_k(x) R_{\nu}(x, hD_x; z, h)$ (see the definition of the global parametrix in Section 3) and $A_{j, k}$ are the operators $A_j$ corresponding to the neighborhood of $x_k.$ It is clear that the integration in $(I)$ and $(II)$ is over compact sets, respectively in $(x, \xi)$ and $(x, \xi')$. On the other hand, an application of Lemma 6 yields
$$(III)= {\mathcal O} (h^{-n + 1-4\epsilon}).$$

Therefore, by using (\ref{eq:4.8}) and (\ref{eq:4.11}) and the action of $(Q - z)^{-1}$, we are going to calculate the trace of
$$ -\frac{1}{\pi}\int_{|\im z | \geq h^{\epsilon}} \one_{\Omega}(x)\partial \tilde{f}(z) R(x, hD_x; z, h) \chi_k(x, \xi') L(dz).$$
The principal symbol of $-\frac{1}{\pi} \int_{|\im z| \geq h^{\epsilon}}\partial\tilde{f}(z) R_{\nu}(x, hD_x; z, h) L(dz)$ is $f(q_{\nu}(x, \xi))$, so the principal symbol of
$$ -\frac{1}{\pi}\int_{|\im z | \geq h^{\epsilon}} \one_{\Omega}(x)\partial\tilde{f}(z) R_{\nu}(x, hD_x; z, h) \chi_k(x, \xi') L(dz)$$
 is $f(q_{\nu}(x, \xi)) \chi_k(x, \xi'),\: x \in \Omega.$ Next the kernel $K_{\nu, k}(x, y)$ of the operator with symbol $f(q_{\nu}(x, \xi)) \chi_k(x, \xi')$ is continuous in $\Omega \times \Omega$ and its trace is 
equal to  $\int_{\Omega} K_{\nu, k} (x, x) dx$. Thus we conclude that

$$ (I) + (II) = (2 \pi h)^{-n} \Bigl(\int_{\Omega} \int f(q_1(x, \xi)) dx d\xi + \int_{\Omega}  \int f(q_2(x, \xi)) dx d\xi\Bigr) + {\mathcal O}(h^{- n + 1}),\: 0 < h \leq h_0(\delta).$$

Consequently,
\begin{eqnarray} \label{eq:4.12}
{\rm tr}(f(Q) )= (2 \pi h)^{-n} \Bigl(\int_{\Omega} \int f(q_1(x, \xi)) dx d\xi + \int_{\Omega}  \int f(q_2(x, \xi)) dx d\xi\Bigr)\nonumber\\
 + {\mathcal O}(h^{- n + 1-4\epsilon}),\: 0 < h \leq h_0(\delta).
\end{eqnarray}

The above argument implies easily that the spectrum of $Q$ in $[0, (\im \om + \delta)^2]$ is formed by isolated eigenvalues of finite multiplicities. To prove this it is sufficient to show that the resolvent $(Q- z_0)^{-1}$ for fixed $z_0$ with $\im z_0 > 0$ is compact. 
Every term $(Q- z_0)^{-1} \chi_k(x, \xi)$ is a trace class operator. On the other hand, for $C_0 > 0$ large enough the ellipticity of boundary problem
and our construction in Section 3 of a (microlocal) parametrix show that $(Q - z_0)^{-1} ( 1 -\chi(h D_{x'}))\Phi_k(x)$ is a compact operator.

\section{Proof of Theorem 1}

In this section we prove the estimate (\ref{eq:1.4}). Let $\re \omega_0 =  1, \: \im \omega_0 > 0.$ We fix $0 < \im \omega_0 < 1$  so that (\ref{eq:cond}) holds and choose $0 <\delta < \im \omega_0$ small enough. Thus the construction of the semi-classical parametrix in Section 3 works and we may exploit the results of Section 4.  Set $(6\delta \im \omega_0 + 9\delta^2) =\delta_1^2$ and 
let $f \in C_0^{\infty}(]-\infty, (\im \omega_0 + 3\delta)^2[ ; [0, 1])$ be a function which is equal to 1 on $[0, (\im \omega_0 + 2\delta)^2]$. Consider $f(q_1(x, \xi))$ and $f(q_2(x, \xi))$, where the symbols $q_j(x, \xi),\: j = 1,2,$ were introduced in Section 3. On the support of $f(q_1(x, \xi))$ we have
$$||\xi|^2 - 1|^2 + (\im \omega_0)^2 \leq (\im \omega_0 + 3\delta)^2$$
which yields
$$||\xi|^2 - 1|\leq ((\im\omega_0) + 3\delta)^2 - (\im \omega_0)^2)^{1/2} = \delta_1.$$
Thus on the support of $f(q_1(x, \xi))$ we have $1 -\delta_1 \leq |\xi|^2 \leq 1 + \delta_1.$
Similarly on the support of $f(q_2(x, \xi))$ we have
$$1 -\delta_1 \leq \fm |\xi|^2 \leq 1 + \delta_1.$$
In the following we suppose that $4 \epsilon < 1/2$ so $1 - n - 4 \epsilon > 1/2 -n.$
Taking the integral over $\Omega \times \R^n$, we deduce
$$\int_{\Omega}\int f(q_1(x, \xi)) dx d\xi \leq  {\rm vol}\:(\Omega) \omega_n\Bigl( (1+ \delta_1)^{n/2} - (1 - \delta_1)^{n/2}\Bigr),$$
 $\omega_n$ being the volume of unit ball in $\R^n$
 and
$$\int_{\Omega} \int f(q_2(x, \xi) )dx d\xi \leq \int_{\Omega} (1 + m(x))^{n/2} dx \omega_n\Bigl( (1+ \delta_1)^{n/2} - (1 - \delta_1)^{n/2}\Bigr).$$

Consequently,
$$\tr f(Q) = (2\pi h)^{-n}\Bigl[ \int_{\Omega}\int f(q_1(x, \xi)) dx d\xi + \int \int f(q_2(x, \xi)) dx d\xi] + {\mathcal O}(h^{1 - n -4 \epsilon})$$
$$= \frac{(1+ \delta_1)^{n/2} - (1 - \delta_1)^{n/2}}{(2 \pi)^n}\omega_n \int_{\Omega} \Bigl( 1 + (1 + m(x))^{n/2}\Bigr) dx   h^{-n} + {\mathcal O}(h^{1 - n -4 \epsilon})$$
$$ = \Bigl((1+ \delta_1)^{n/2} - (1 - \delta_1)^{n/2}\Bigr)C_2(\Omega, n) h^{-n} + {\mathcal O}(h^{1-n - 4\epsilon}) : = {\mathcal N}(2\delta; h).$$
with $C_2(\Omega, n) = \frac{\omega_n}{(2 \pi)^n} \int_{\Omega} \Bigl( 1 + (1 + m(x))^{n/2}\Bigr) dx.$
 This implies that the number $M(2\delta; h)$ of the eigenvalues of $|Q|^{1/2}$ smaller than $\im \omega_0 + 2\delta$ satisfies
\begin{equation}
M(2\delta; h) \leq {\mathcal N}(2\delta; h). 
\end{equation}

Let $0 \leq \mu_1 \leq \mu_2 \leq ...$ be the eigenvalues of $|Q|^{1/2}$ followed by an infinite repetition of $\inf\sigma(|Q|^{1/2})$ in the case if there are only  finitely many eigenvalues. Since $Q- z$ is related to an elliptic boundary problem for $z < (\im \omega_0)^2$, we have
$$\mu_1 \geq \im \omega_0 - o(1), \: h \to 0.$$
Now let $N = N(\delta; h)$ be the number of the eigenvalues $\nu_1, \nu_2,...,\nu_N$ of $h^2{\mathcal P} - \omega_0$ with $0 \leq |\nu_1|\leq ...\leq |\nu_N|\leq \im \omega_0 +\delta.$ \\

Applying the Weyl inequality for the non-selfadjoint operator $h^2{\mathcal P} - \omega_0$ (see Appendix a. in \cite{Sj}), we obtain
\begin{equation}
\mu_1 ... \mu_N \leq |\nu_1|...|\nu_N|.
\end{equation}
As in \cite{SjZ}, this leads to the estimate
\begin{equation}
N(\delta; h) \leq c_0 M(2\delta; h) \leq c_0 {\mathcal N}(2\delta; h), \: 0 < h \leq h_0(\delta)
\end{equation}
with a constant $c_0 = c_0(\omega_0, \delta)> 0$ depending only on  $\omega_0$ and $\delta.$ For reader convenience, we present the argument. Assume that
$M(2\delta; h) \leq N(\delta; h)$. Then
$$\mu_1^M (\im \omega_0 + 2\delta)^{N - M} \leq (\im \omega_0 + \delta)^N$$
and we get
$$\Bigl[ \frac{\im \omega_0 + 2\delta}{ \im \omega_0 + \delta} \Bigr]^N \leq \Bigl[\frac{\im \omega_0 + 2\delta} {\mu_1} \Bigr]^{M}.$$
Therefore for $h \leq h_0(\delta)$ we have
$$N \log \Bigl(\frac{\im \omega_0 + 2\delta}{ \im \omega_0 + \delta} \Bigr) \leq M \log \Bigl( \frac{\im \omega_0 + 2\delta} {\im \omega_0 - \delta}\Bigr).$$
We may estimate
$$ c_0 = \frac{ \log \Bigl( 1 + \frac{3 \delta}{\im \omega_0 - \delta}\Bigr)}{\log \Bigl( 1 + \frac{\delta}{\im \omega_0 + \delta}\Bigr)}$$
as $\delta \to 0$ and conclude that with $c_0 = 3(1 + {\mathcal O}_{\omega_0}(\delta))$ we have
 $$N(\delta; h) \leq c_0 M( 2\delta; h).$$

 Consider the disk $D(\omega_0; \im \omega_0 + \delta)$ with center $\omega_0 \in \C$ and radius $\im \omega_0 + \delta < 1.$ 
Let $S_0$ be the rectangle
$$S_0 = \{(x, y)\in \R^2:\: 1 -\delta_2 \leq x \leq 1 + \delta_2,\: -X \leq y \leq \im \omega_0 \} \subset D(\omega_0; \im \omega_0 + \delta),$$
where
$$\delta_2 = \sqrt{1 - \delta^2}\Bigl(2 \im \omega_0 \delta + \delta^2\Bigr)^{1/2} < \delta_1,$$
and  $X = \Bigl((\im \om)^2 + \delta^2(2 \im \om \delta + \delta^2)\Bigr)^{1/2} - \im \om > 0.$ For $\delta$ small enough we have $X < \im \om.$ We choose $0 < \theta < \pi/ 2$ so that 
$$\tan \theta = \frac{X}{1 + \delta_2}  < X.$$
It is clear that when $\theta \searrow 0$ we have $\delta \searrow 0.$
According to the above result, the number of the eigenvalues of the operator $h^2 {\mathcal P}$ in $S_0$ can be estimated by
$c_0{\mathcal N}(2\delta; h),$ so the number of the eigenvalues of ${\mathcal P}$ in $h^{-2} S_0$ is estimated by $c_0{\mathcal N}(2\delta; h).$

Given $r \gg 1$, consider the set
$$S_{0, K} =\bigcup_{k = 0}^K \Bigl(\frac{1 + \delta_2}{1 - \delta_2}\Bigr)^{-k} r S_0,$$
where $ K = K(r)$ is chosen so that $\Bigl(\frac{1+ \delta_2}{1- \delta_2}\Bigr)^{K} r^{-1} \leq h_0^2(\delta)$. Therefore for
$$h = \Bigl(\frac{1  + \delta_2}{1 - \delta_2}\Bigr)^{k/2} r^{-1/2} \leq h_0(\delta)$$ 
we can apply the above estimate and the number of the eigenvalues of ${\mathcal P}$ in $S_{0, K}$ is estimated by
$$ c_0\sum_{k = 0}^K {\mathcal N}\Bigl(2\delta; \Bigl(\frac{1 + \delta_2}{1 - \delta_2}\Bigr)^{k/2} r^{-1/2}\Bigr)\leq \Bigl((1+ \delta_1)^{n/2} - (1 - \delta_1)^{n/2}\Bigr) c_0C_2(\Omega, n) \sum_{k= 0}^K \Bigl(\frac{1 - \delta_2}{1 + \delta_2}\Bigr)^{nk/2} r^{n/2} $$
$$+  \sum_{k=0}^{K} \Bigl(\frac{1 - \delta_2}{1 + \delta_2}\Bigr)^{(n- 1 +4\epsilon)k/2}{\mathcal O} (r^{(n- 1 + 4\epsilon)/2})$$
$$\leq c_0C_2(\Omega, n) (1+ \delta_2)^{n/2}\frac{(1+ \delta_1)^{n/2} - (1 - \delta_1)^{n/2}}{(1 + \delta_2)^{n/2} - (1 - \delta_2)^{n/2}} r^{n/2} + {\mathcal O}_{\delta}(r^{(n- 1+4 \epsilon)/2}), \: r \geq r_0( \delta),$$

Clearly,
$$\Bigl((1 + \delta_1)^{n/2} - (1 - \delta_1)^{n/2}\Bigr) =n \delta_1 ( 1 + {\mathcal O}_{\om}(\delta)),$$
$$\Bigl((1 + \delta_2)^{n/2} - (1 - \delta_2)^{n/2}\Bigr) = n\delta_2 ( 1 + {\mathcal O}_{\om}(\delta))$$
and 
$$\frac{\delta_1}{\delta_2} = \frac{1}{\sqrt{1 - \delta^2}}\Bigl(\frac{6 \im \om + 9 \delta}{ 2 \im \om + \delta}\Bigr)^{1/2} \to_{\delta \to 0} \sqrt{3}.$$
 Thus the number of the eigenvalues of ${\mathcal P}$ in $S_{0, K}$ can be estimated by
$$3\sqrt{3} (1 + {\mathcal O}_{\omega_0}(\delta)) \frac{\omega_n}{(2\pi)^n} \int_{\Omega} \Bigl( 1 + (1 + m(x))^{n/2}\Bigr)dx r^{n/2}, \: r \geq r_0(\delta),$$
where we increase $r_0(\delta)$, if it is necessary.
On the other hand , it is easy to see that for small $\theta > 0$ we have
$$\Lambda_{\theta, r} = \{ z \in \C:\: |\arg z |\leq \theta,\: |z| \leq r\} \subset \bigcup_{k= 0}^{\infty}\Bigl(\frac{1 + \delta_2}{1 - \delta_2}\Bigr)^{-k} r S_0$$ 
$$= S_{0, K} \cup \bigcup_{k= K+1}^{\infty}\Bigl(\frac{1 + \delta_2}{1 - \delta_2}\Bigr)^{-k} r S_0 = S_{0, K} \cup Q_{K}.$$
Here $Q_K$ is a compact set independent on $r$. Indeed, for $j \geq  K(r) + 1$ we have 
$$\Bigl(\frac{1 + \delta_2}{1 -\delta_2}\Bigr)^{-j} r < \frac{1}{h_0^2(\delta)},$$
hence the rectangles $\Bigl(\frac{1 + \delta_2}{1 - \delta_2}\Bigr)^{-j} r S_0$ for $ j \geq K(r) + 1$ are included in a fixed compact set.

 Finally, the number $N(\theta, r)$ of the eigenvalues of ${\mathcal P}$ in $\Lambda_{\theta, r}$ is estimated by
\begin{eqnarray} 
N(\theta, r) \leq 3\sqrt{3}( 1 + {\mathcal O}_{\omega_0}(\delta)) \frac{\omega_n}{(2\pi)^n}\Bigl(\int_{\Omega} \Bigl( 1 + (1 + m(x))^{n/2}\Bigr)dx \Bigr) r^{n/2} + C_0(\delta), \:r \geq r_0( \delta).
\end{eqnarray}
As we mentioned above, as $\theta \to 0$, we have $\delta \to 0$. Hence for fixed $\omega_0$ we can replace $\Oc_{\omega_0}(\delta)$ by $\ep(\theta) \to 0$ as $\theta \to 0$ and $r_0(\delta)$ by $r(\theta).$
This completes the proof of Theorem 1.

{\footnotesize

\end{document}